\documentclass[12pt,final,a4paper]{article}
\pdfoutput=1
\usepackage{amsfonts}
\usepackage{amsmath,amssymb,latexsym}
\usepackage{mathtools}
\DeclarePairedDelimiter{\opnorm}{\lVert}{\rVert_{1}}
\DeclarePairedDelimiter{\norm}{\lVert}{\rVert}
\usepackage{ulem}
\usepackage[
tracking=true, %
kerning=true, %
spacing=true, %
factor=1100, %
stretch=10, %
shrink=10 %
]{microtype}
\usepackage{epsf}
\usepackage{color}
\usepackage{csquotes}
\usepackage{epsfig}
\usepackage{graphicx}
\usepackage{mathdots}
\usepackage{array}

\usepackage{subcaption}
\usepackage{tikz,pgfplots}
\usetikzlibrary{arrows.meta}
\usetikzlibrary{backgrounds}
\usepgfplotslibrary{patchplots}
\usepgfplotslibrary{fillbetween}
\pgfplotsset{%
	width=\textwidth,
	height=0.28\textheight,
	tick label style={font=\scriptsize},
	title style={font=\scriptsize},
	legend style={font=\tiny, at={(0.0,1.05)}, anchor=south west},
	legend columns=10,
	xlabel shift={-1.5ex},
	xlabel style={font=\footnotesize},
	ylabel style={font=\scriptsize},
	ylabel shift={-0.75em},
	layers/standard/.define layer set={%
		background,axis background,axis grid,axis ticks,axis lines,axis tick labels,pre main,main,axis descriptions,axis foreground%
	}{
		grid style={/pgfplots/on layer=axis grid},%
		tick style={/pgfplots/on layer=axis ticks},%
		axis line style={/pgfplots/on layer=axis lines},%
		label style={/pgfplots/on layer=axis descriptions},%
		legend style={/pgfplots/on layer=axis descriptions},%
		title style={/pgfplots/on layer=axis descriptions},%
		colorbar style={/pgfplots/on layer=axis descriptions},%
		ticklabel style={/pgfplots/on layer=axis tick labels},%
		axis background@ style={/pgfplots/on layer=axis background},%
		3d box foreground style={/pgfplots/on layer=axis foreground},%
	},
}
\tikzset{
	every node/.append style={
		font=\footnotesize %
	},
	dashed/.style={dash pattern=on 4pt off 4pt}
}

\makeatletter
\def\Ddots{\mathinner{\mkern1mu\raise\p@
		\vbox{\kern7\p@\hbox{.}}\mkern2mu
		\raise3\p@\hbox{.}\mkern2mu\raise5\p@\hbox{.}\mkern1mu}}
\makeatother
\usepackage{tabularray}
\SetTblrInner{rowsep=0.1ex}
\UseTblrLibrary{siunitx,booktabs,counter}
\newcounter{mycnta}

\usepackage{booktabs}
\usepackage{caption}
\usepackage{multirow}
\usepackage{siunitx}
\usepackage{hyperref}

\newcommand{\blue}[1]{#1}%

\def\en{\text{\sc{e}-}}

\def\mdim{N} %

\def\id{I}

\newcommand{\e}{\mathrm{e}}

\def\beq{\begin{equation}}
\def\eeq{\end{equation}}

\begin{document}
\title{Efficient scaling and squaring method for the matrix exponential}

\author{Sergio Blanes\thanks{Instituto de Matem\'atica Multidisciplinar,
		Universitat Polit\`ecnica de Val\`{e}ncia, E-46022  Valencia, Spain. Email:  \texttt{serblaza@imm.upv.es}}
	\and
	Nikita Kopylov\thanks{Instituto de Matem\'atica Multidisciplinar,
		Universitat Polit\`ecnica de Val\`{e}ncia, E-46022  Valencia, Spain. Email:  \texttt{nikop1@upv.es}}
	\and
	Muaz Seydao\u{g}lu\thanks{Department of Mathematics, Faculty of Art and Science, Mu\c{s} Alparslan University, 49100, Mu\c{s}, Turkey. Email: \texttt{m.seydaoglu@alparslan.edu.tr}}
}

\maketitle

\begin{abstract}
	This work presents a new algorithm to compute the matrix exponential within a given tolerance.
	Combined with the scaling and squaring procedure, the algorithm incorporates Taylor, partitioned and classical Padé methods shown to be superior in performance to the approximants used in state-of-the-art software.
	The algorithm computes matrix--matrix products and also matrix inverses, but it can be implemented to avoid the computation of inverses, making it convenient for some problems.
	If the matrix \( A \) belongs to a Lie algebra, then \( \e^A \) belongs to its associated Lie group, being a property which is preserved by diagonal Pad\'e approximants, and the algorithm has another option to use only these.
	Numerical experiments show the superior performance with respect to state-of-the-art implementations.


\end{abstract}

\section{Introduction}
\label{sec.1}
We present an algorithm for computing the matrix exponential $\e^A$ \cite{moler03ndw} within a specified tolerance $tol$, given a complex matrix $A$ of size $\mdim\times\mdim$.
First, the algorithm computes a bound, $\theta$, to the norm of the matrix, i.e. $\opnorm{A}<\theta$ (the 1-norm can be replaced by any other norm), and then, according to \( \theta \), it chooses the scheme among a list of selected methods which provides an approximation to $\e^A$ with in such tolerance.
If no matching scheme is found, the algorithm employs scaling and squaring (summarized further) to find the cheapest method matching $\tfrac{\opnorm{A}}{2^{s}}<\theta$.
It is observed that different schemes for approximating the matrix exponential should be selected based on the specific problem at hand.

As an example of numerical methods where such adaptivity is beneficial, consider exponential integrators for solving differential equations.
They have been shown to be very useful for a significant number of problems \cite{hochbruck10ein}.
For instance, most Lie-group methods (see also \cite{iserles00lgm} for a review) like Magnus integrators (see \cite{blanes09tme} and references therein), Crouch--Grossman methods \cite{crouch93nio}, Runge--Kutta--Munthe-Kaas methods \cite{munthekaas98rkm}, etc. require calculating of one or several matrix exponentials per step.
In most cases, it suffices to approximate each matrix exponential only up to a given tolerance lower than the round-off accuracy, which would make integration computationally more expensive.
When the Lie-group structure must be preserved, one can still provide cheaper approximations to the exponential still preserving the structure.

For example, let us consider the nonlinear differential equation
\begin{equation}\label{eq.Crouch-Grossman0}
	Y' = A(t,Y) \, Y, \qquad Y(0)=Y_0 \in G,
\end{equation}
where $A\in\mathbb{R}^{N\times N}$ and $G$ is a matrix Lie group,
which appears in relevant physical fields such as rigid
body mechanics, in the calculation of Lyapunov exponents and other problems arising in Hamiltonian dynamics.
It can be shown that a differential equation evolving on a matrix Lie group $G$ can be written in the form \eqref{eq.Crouch-Grossman0} \cite{iserles00lgm}.
In \cite{crouch93nio} the authors present, e.g. the following 3-stage method of order 3
\begin{equation} \label{eq.Crouch-Grossman1}
	\begin{aligned}
		Y^{(1)} & = Y_n,                                                     &  & K_1=A(t_n,Y^{(1)}),      \\
		Y^{(2)} & = \mathrm{e}^{ha_{2,1}K_{1}} \, Y_n,                      &  & K_2=A(t_n+c_2h,Y^{(2)}), \\
		Y^{(3)} & = \mathrm{e}^{ha_{3,2}K_{2}} \, \mathrm{e}^{ha_{3,1}K_{1}} \, Y_n, &  & K_3=A(t_n+c_3h,Y^{(3)}), \\
		Y_{n+1} & = \mathrm{e}^{hb_3 K_{3}}\,\mathrm{e}^{hb_2 K_{2}} \,\mathrm{e}^{hb_1 K_{1}} \, Y_n, &  &
	\end{aligned}
\end{equation}
with
\begin{gather*}
	a_{2,1}=\frac34, \quad
	a_{3,2}=\frac{17}{108}, \quad
	a_{3,1}=\frac{119}{216}, \quad
	b_{1}=\frac{24}{178}, \quad
	b_{2}=-\frac23, \quad
	b_{3}=\frac{13}{51},\\
	c_2=a_{2,1},\quad
	c_3=a_{3,1}+a_{3,2},
\end{gather*}
which requires the computation of six matrix exponentials per step.
If the norm of $A$ does not change much during each step we will have that $\opnorm{K_{1}}\sim\opnorm{K_{2}}\sim\opnorm{K_{3}}$ on each time step, but the coefficients $a_{i,j}$ and $b_i$ differ up to nearly one order of magnitude.
For this problems it is then obvious that, since the method is only of third order: {\it (i)} it is not necessary to approximate the matrix exponentials up to round off accuracy as far as the approximations retain the desired structure, and {\it (ii)} each matrix exponential can be approximated by a different scheme to reach the desired tolerance (which can also change at different steps) in order to reduce the computational cost while retaining the accuracy and qualitative properties.

As a second example, consider methods for stochastic differential equations (SDEs), such as \cite{Ta2015} that computes matrix exponentials explicitly.
For this type of equations, integrators are usually of low order and involve repeating the simulation a sufficiently large number of times.
Therefore, similarly to the example above, it is beneficial to use cheaper schemes to lower the overall simulation time.

On the other hand, in \cite{lezcanocasado20aam} it is claimed that matrix multiplications without inverse make the exponential of matrices a very cheap function to use in the context of deep learning, and a Taylor algorithm has been added to PyTorch 1.7.0 showing a significant speed improvement over Tensorflow.
For those problems, polynomial approximations to the matrix exponential are advantageous.

Scaling and squaring is perhaps one of the most frequently used approaches in popular computing packages such as \textsc{MATLAB}'s \texttt{expm}, Mathematica's \texttt{MatrixExp}, and Julia's \texttt{exp} and \texttt{ExponentialUtilities.jl} in combination with Pad\'e approximants \cite{almohy09ans,higham05tsa,higham10cma,julia17,Rackauckas2017}.
Specifically, it is based on the key property
\begin{equation}  \label{eses}
	\e^A = \left( \e^{A/2^s}  \right)^{2^s}
	, \quad s\in\mathbb{N}.
\end{equation}
The exponential $\e^{A/2^s}$ is replaced by an approximation $w(A/2^s)$, which usually is chosen to approximate the exponential up to round off error, and it is then squared $s$ times.

The goal of this work is to construct an algorithm that, when provided with a matrix \( A \in \mathbb{C}^{\mdim \times \mdim} \) and a tolerance \( tol \), computes a function \( w_{\alpha}(A) \), where $\alpha$ refers to a label to identify the method, capable of approximating the matrix exponential \( \e^A \) so that
\begin{equation}\label{eq.relerr}
	relerr\coloneq\frac{\opnorm{w_{\alpha}(A)-\e^A}}{\opnorm{\e^{A}}} < tol \cdot \opnorm{{A}}.
\end{equation}

The most popular choices for the function $w_{\alpha}(A)$ are:
\begin{itemize}
	\item Taylor polynomial: $w_{\alpha}(A)=t_m(A)$, where $t_m(x)=\sum_{n=0}^{m}\tfrac{x^n}{n!}, \ x\in \mathbb{C}$, and $|t_m(x)-\e^x|=\mathcal{O}(x^{m+1})$.
	\item Rational Pad\'e approximant:  $w_{\alpha}(A)=r_{k,m}(A)$, where $r_{k,m}(x)=p_{k,m}(x)/q_{k,m}(x)$ with
	\begin{equation}  \label{eq.1.2b}
		p_{k,m}(x)= \sum_{j=0}^k \frac{(k+m - j)! k!}{(k+m)! (k-j)!} \frac{x^j}{j!},
		\quad
		q_{k,m}(x)= \sum_{j=0}^m \frac{(k+m - j)! m!}{(k+m)! (m-j)!} \frac{(-x)^j}{j!},
	\end{equation}
	which is an approximation of order $s=k+m$ with the leading error term given \cite{Iserles1982, higham09tsa} by
	\begin{equation}\label{eq.ErrPade}
		\e^x - r_{k,m}(x)=(-1)^m \frac{ k! m!}{(k+m)! (k+m+1)!}x^{k+m+1} + \mathcal{O}(x^{k+m+2}).
	\end{equation}
	Notice that the Taylor polynomial is a special case of Padé approximant,  $t_s(x)=r_{s,0}(x)$, and has a larger leading error than $r_{k,m}(x)$ with $s=k+m$, $k,m>0$.
	However, the performance of a method strongly depends on the computational cost\footnote{In some cases round off accuracy must be taken into account.} and it should also be considered.

	\item Rational Chebyshev approximants:  $w_{\alpha}(A)=\tilde r_{k,m}(A)$, where $\tilde r_{k,m}(x)=\tilde p_{k}(x)/\tilde q_{m}(x)$, and $\tilde p_{k}(x),\ \tilde q_{m}(x)$ are polynomials of degree $k$ and $m$, respectively, with $k\leq m$ and with coefficients chosen such that
	\[
	|\tilde c_{k,m}(x)-\e^{x}|< tol_{k,m}, \quad x\in(-\infty,0) .
	\]
	where $tol_{k,m}$ depends on the values of $k$ and $m$, and the choice of the coefficients for the polynomials.

	These methods are addressed for the case in which the matrix $A$ is diagonalizable with negative real eigenvalues with some of them taking very large values, and has received much interest \cite{cody69cra,trefethen06tqa}. This problem does not fit in the class of problems we consider in this work ($\|A\|$ must be of moderate size), but most techniques used in this work apply, and this problem will be considered in the future.
\end{itemize}
We analyse an extensive set of Taylor and rational Pad\'e methods, we get error bounds for a set of tolerances and, taking into account their computational cost, we select a list of methods which our algorithm uses according to the provided matrix $A$ and the tolerance $tol$.

Among all possible choices, one should use the method that provides the desired accuracy at the lowest computational cost, and this requires to carry two types of analysis:
\begin{description}
	\item[Backward error.]
	Given a particular method $w_{\alpha}$, we look for an associated scalar function for that method, say $\theta_{\alpha}(y)$, such that, given a matrix $A$ and a positive integer number, $s$, such that if $\opnorm{A} \leq 2^{s}\theta_{\alpha}(tol)$, the method provides an approximation with a relative error below the tolerance when used with $s$ squarings.
	Forward and backward error analysis are frequently used in the literature, and we will consider the backward error analysis in this work, similarly to the analysis carried out in \cite{higham05tsa,higham10cma} and implemented in the function (\texttt{expm}) in \textsc{MATLAB}.

	\item[Computational cost.]
	For the analysis of the cost we assume that one dense matrix--matrix multiplication cost is $C$, which we will take as the reference cost.
	We will say that the cost of a method is $k \, C$ if its computational cost is approximately $k$ times the cost of a matrix--matrix product.
	Then, given two matrices $A$ and $B$, the cost to compute $A^{-1}$ or $A^{-1} B$ will be taken as $\tfrac{4}{3} C$.
	For a shorter notation, we assume $C=1$.

	Contrarily to the error analysis, the study of computational cost can be considered as an art rather than a proper theory.
	During the last decades, the best existing method to solve a particular problem has been changing as new techniques to reduce the cost of the algorithms are found (and will continue in the future with the new algorithms and computer architectures).
	For example, Taylor methods were originally discarded since, when computed with the standard Horner's algorithm, they require $m-1$ products to compute $t_m$.
	However, $t_m$ with $m=2,4,6,9,16,20$ can be computed with $k=1,2,3,4,5,6$ products, respectively, a significant saving which made the Taylor methods competitive \cite{paterson73otn,sastre15nsq}.
	In addition, in \cite{preprint,bader19ctm,sastre18eeo} it is shown that further reduction can be carried out such that the Taylor polynomial with $m=2,4,8,12,18$ can be computed with $k=1,2,3,4,5$ products, respectively, making them the methods of choice for many problems.

	On the other hand, diagonal Pad\'e methods have the property that $q_{m,m}(x)=p_{m,m}(-x)$, and this symmetry allows finding a procedure to compute both polynomials $p_{m,m}(x)$ and $ q_{m,m}(x)$ simultaneously at a reduced cost for $m=1,2,3,5,7,9,13$ with $k=0,1,2,3,4,5,6$ products, respectively, in addition to one inverse to compute $r_{m,m}$.
	In \cite{higham05tsa} it is claimed that, since $r_{m,m}$ can be computed with the same cost as $r_{k,m}$ or $r_{m,k}$ with $k<m$ and provide higher accuracy (see eq. \eqref{eq.ErrPade}), only diagonal Pad\'e methods are considered.
	However, in \cite{sastre12emr} it is shown that this is not necessarily the case, showing that a number of approximants $r_{k,m}$, with $k>m$ and with an appropriate fractional decomposition, can be computed at the same cost as $r_{m,m}$ while providing higher accuracy.
\end{description}

In addition, when analyzing the methods, the following restrictions should also be taken into account:
\begin{description}
	\item[Avoiding complex coefficients.]
	In this work we will decompose some rational Pad\'e methods into simpler fractions, and some caution must be taken to avoid schemes with complex coefficients.
	The goal of this work is to propose methods which are optimized for real matrices and also working efficiently for complex matrices, so we look for schemes with real coefficients.
	Note that the list of methods employed in our algorithm can be easily extended by adding more methods with complex coefficients addressed for complex matrices.
	At this point we will take into account that for $m=2j$ the denominators in $r_{k,m}$, say $q_{k,2j}, \ j=1,2,\ldots$ in the Pad\'e approximation to the matrix exponential have $j$ distinct complex roots, $\alpha_i$, and their conjugates,  $\bar{\alpha}_i$, so
	\[
	q_{k,2j}=\delta \prod_{i=1}^{j}q^{(k)}_{j,i}, \qquad q^{(k)}_{j,i}=(x-\alpha_i)(x-\bar\alpha_i)=(x^2-2Re(\alpha_i)+|\alpha_i|^2),
	\]
	and $\delta=(-1)^mk!/(k+m)!$ In the case that  $m=2j+1, \ j=0,1,2,\ldots$, we have that $q_{k,2j+1}$ has only one real root and $j$ complex ones with their conjugates. This property will be taken into account when decomposing $r_{k,m}$ into simpler fractions.

	\item[Lie group methods.]
	If $A$ belongs to a Lie algebra, then $\e^A$ will belong to the associated Lie group.
	If this property should be preserved in a problem, the presented algorithm allows computing the exponential using only diagonal Pad\'e approximants which preserve this property to round off accuracy.
\end{description}

Eventually, the procedure to find the most efficient schemes is as follows:
\begin{itemize}
	\item
	To lower the computational cost, some Pad\'e approximants are decomposed into simpler fractions avoiding schemes with complex coefficients.
	\item
	For each selected Taylor and Pad\'e method, say $w_{\alpha}(x)$, we numerically obtain the function $\theta_{\alpha}(t)$ such that if \( \opnorm{A}\leq\theta_{\alpha}(tol) \), then \autoref{eq.relerr} is satisfied for $tol=10^{-k}, \ k=0,1,\ldots,16$.
	If this condition is not satisfied, the scaling and squaring procedure is applied.
\end{itemize}

As a result of analysis, in \autoref{sec.list} we will provide a list of Taylor and Pad\'e methods, each one being the optimal one for at least some values of the tolerance and norm of the matrix.

\section{Backward error analysis}
In the implementation of the scaling and squaring algorithm, the choice of the optimal method and the scaling parameter for a given matrix $A$ and tolerance, $tol$, are based on the control of the backward error \cite{higham05tsa}.
More specifically, given an approximation $w_n(A)$ of order $n$ to the exponential , i.e.
$w_n(x) = \e^x + \mathcal{O}(x^{n+1})$, one defines the function $h_{n+1}(x) = \log(\e^{-x} w_n(x))$, then
$w_n(2^{-s} A) = \e^{2^{-s} A + h_{n+1}(2^{-s} A)}$ and
\[
\left( w_n(2^{-s}A) \right)^{2^s} = \e^{A + 2^s h_{n+1}(2^{-s} A)} \equiv \e^{A + \Delta A},
\]
where $\Delta A \coloneq 2^s h_{n+1}(2^{-s}A) $ is the backward error originating in the approximation of $\e^A$.
If in addition $h_{n+1}$ has a power series expansion
\[
h_{n+1}(x) = \sum_{k=n+1}^{\infty} c_k x^k,
\]
with a non-zero radius of convergence, then it is clear that $\|h_{n+1}(A)\| \le \widetilde{h}_{n+1}(\|A\|)$, where
\[
\widetilde{h}_{n+1}(x) = \sum_{k=n+1}^{\infty} |c_k| x^k,
\]
and thus, taking into account the scaling of \( A \), one gets
\begin{equation}  \label{bound.1}
	\frac{\|\Delta A \|}{\| A\|}=\frac{\|h_{n+1}(2^{-s}A)\|}{\|2^{-s}A\|} \le \frac{\widetilde{h}_{n+1}(\|2^{-s}A\|)}{\|2^{-s}A\|}.
\end{equation}
Given a prescribed tolerance $tol$, one computes numerically
\begin{equation}  \label{bound.2}
	\theta_n(tol) = \max \left\lbrace  \theta : \frac{\widetilde{h}_{n+1}(\theta)}{\theta} \le tol \right\rbrace,
\end{equation}
for different values of $tol$ computed by truncating the series of the corresponding functions
$\widetilde{h}_{m+1}(\theta)$ after 150 terms.
Then $s$ is chosen so that $\| 2^{-s} A\| \le \theta_n(tol)$ and $\left( w_n(2^{-s} A) \right)^{2^s}$ is used to approximate $\e^A$.
Notice that for small values of $\| \Delta A\|$ we have that
\[
\norm{\Delta A} = \frac{\norm{w_n(A)-\e^A}}{\norm{\e^A}} + \mathcal{O}( \norm{\Delta A}^2).
\]

In the particular case that $w_n$ is a diagonal Pad\'e approximant $r_{m,m}$, then $n = 2m$.
The values of $\theta_{2m}(tol)$ for $tol=u$ where $u$ is unit round-off in single and double precision are shown in \cite{almohy09ans} and collected in \autoref{tab.theta} for convenience of the reader.
It should be noted that \texttt{expm} in \textsc{MATLAB} is implemented only for the case $u=2^{-53}\simeq 1.1\cdot 10^{-16}$.
According to Higham  \cite{higham05tsa}, $m=13$ and therefore $r_{13,13}$ is the optimal choice in double precision when scaling is required.
When $\|A\| \le \theta_{26}(2^{-53})=5.37$, the algorithm in  \cite{higham05tsa} takes the first $m \in \{3,5,7,9,13\}$ such that $\|A\| \le \theta_{2m}(u)$.
We have extended the analysis to other rational approximations and tolerances.
\autoref{tab.thetaAll} shows the results obtained only for a selected set of values of the tolerance for simplicity in the presentation.
Only those methods which show the best performance for some values of $\norm{A}$ and $tol$ are considered.
Diagonal Pad\'e methods are listed in a separate \autoref{tab.thetaDiag}.

\begin{table}\centering\footnotesize
	\caption{\label{tab.theta} Values of $\theta_{2m}$ for the diagonal Pad\'e approximant $r_m$ of order $2m$ with the minimum number of products $\pi$ for single and double precision. In bold, we have highlighted the asymptotically optimal order at which it is advantageous to apply scaling and squaring since the increase of $\theta$ per extra product is smaller than the factor 2 from squaring.
	}
	\newcolumntype{H}{@{}>{\lrbox0}l<{\endlrbox}}
	\newcolumntype{D}{>{$}r<{$}}
	\begin{tabular}{DDDDHDHDHDHHHD}
		\toprule
		\pi:          & 0        & 1        & 2        & 3        & 3    & 4    & 4    & 5    & 5    & 5             & 6  & 6 & 6 \\
		2m:           & 2        & 4
		& 6        & 8        & 10       & 12       & 14   & 16   & 18   & 20   & 22   & 24            & 26         %
		\\ \midrule
		u \leq2^{-24} &
		8.46\en4      &
		8.09\en2      &
		4.26\en1      &
		1.05          &
		1.88          &
		2.85          &
		\mathbf{3.93} &
		5.06          &
		6.25          &
		7.47          &
		8.71          &
		9.97          &
		11.2                                                                                                                      \\
		u\leq2^{-53}  &
		3.65\en8      & 5.32\en4 &
		1.50\en2      & 8.54\en2 & 2.54\en1 & 5.41\en1 & 9.50\en1 & 1.47 & 2.10 & 2.81 & 3.60 & 4.46 & \mathbf{5.37}              %
		\\
		\bottomrule
	\end{tabular}
\end{table}

Pad\'e approximants are not, of course, the only option one has in this setting.
Another approach to the problem consists in using the  Taylor polynomial of degree $m$  as the underlying approximation $t_m$ to the exponential, i.e.,
taking $ t_m(A) = \sum_{k=0}^{m} A^{k}/k!$ computed efficiently which make them competitive in many cases.

\section{Computational cost}

\subsection{Efficient computation of rational Pad\'e approximants}
\label{sec.4.1}

Diagonal Pad\'e approximants have the form
$ r_{m,m} = \left[ p_{m}(-A) \right]^{-1} p_{m}(A) $ and the evaluation of both $p_{m}(A )$ and $p_{m}(-A )$ can be carried out trying to minimize the number of matrix products.
For illustration, following \cite{higham05tsa,higham09tsa} we reproduce the computation of $r_{5,5}(A)$:
\begin{equation}\label{eq.pade10}
\begin{aligned}
	& u_{5}  =  A \left( b_5 A_4 +  b_3 A_2 + b_1 \id \right), \\
	& v_{5} =  b_4 A_4 + b_2 A_2 + b_0 \id,                 \\
	& (-u_{5}+v_{5}) \, r_{5,5}(A) = u_{5}+v_{5},
\end{aligned}
\end{equation}
with coefficients $b_j$ from \eqref{eq.1.2b}, whereas for $r_{13,13}(A)$ one has
\begin{equation}\label{eq.pade26}
\begin{aligned}
	& u_{13}  = A \left( A_6 \left( b_{13}A_6 + b_{11}A_4 + b_9A_2 \right) + b_7A_6 + b_5A_4 + b_3A_2 + b_1 \id \right), \\
	& v_{13} =  A_6 \left( b_{12}A_6 + b_{10}A_4 + b_8A_2 \right)  + b_6A_6 + b_4A_4 + b_2A_2 + b_0\id,                \\
	& (-u_{13}+v_{13}) r_{13,13}(A) = u_{13}+v_{13}.
\end{aligned}
\end{equation}
Here $A_2=A^2$, $A_4=A_2^2$ and $A_6=A_2A_4$.
Written in this form, it is clear that only three and six matrix multiplications and one inverse are required to obtain approximations of order 10 and 26
to the exponential, respectively.
Diagonal Pad\'e approximants $r_{m,m}(A)$ with $m=3,5,7,9$ and $13$ are used in fact by the function  \texttt{expm} in \textsc{MATLAB}.

However, in many cases it is possible to compute higher order superdiagonal approximants (with $k>m$) at the same computational cost as the diagonal one leading to more efficient schemes \cite{sastre12emr}.
This can be illustrated with the following example. Let
\[
r_{1,1}(x)=\frac{1+\frac12 x}{1-\frac12x},
\]
which can be computed with one inverse.
However, $r_{2,1}(x)$ when decomposed as follows
\[
r_{2,1}(x)=
\frac{1+\frac23x+\frac16x^2}{1-\frac13x} =
 \alpha-\frac12x + \frac{\beta + \gamma x}{1-\frac13x} \eqcolon
 p_0(x) + \frac{p_1(x)}{p_2(x)},
\]
has the same computational cost while having a higher order and resulting in smaller errors for all tolerances of interest.
Since $r_{2,1}(0)=1$ then $ \alpha+\beta=1$ and there is one free parameter to choose.
In general, when the free term of the polynomial is put into the rational function, that is $\alpha=0$ (or, equivalently, $p_0(0)=0$) the round-off errors are lower.
For this reason, in all cases we will consider $p_0(0)=0$.

Our goal is to compute the exponential of matrices with small to moderate norm.
Therefore, using superdiagonal Pad\'e approximants is a valid option.
Similarly, it is easy to see that one can write $r_{2m,m}(x)$ as
\[
r_{2m,m}(x)=p_0(x) + \frac{p_1(x)}{p_2(x)},
\]
$p_i(x), \ i=0,1,2$ being polynomials of degree $m$, which can be computed with $m-1$ products, and one inverse.
Then, at the same cost as the previous $r_{2,2}(x),r_{3,3}(x)$ and $r_{5,5}(x)$ (1, 2, and 3 products, respectively, and one inverse) it is possible to compute $r_{4,2}(x)$, $r_{6,3}(x)$ and $r_{8,4}(x)$, which are more accurate than the previous methods of the same cost for all tolerances of practical interest.

This procedure can be extended to the case $r_{3m,2m}(x)$
\[
r_{3m,2m}(x)=p_0(x) + \frac{p_1(x)}{p_2(x)}+ \frac{p_3(x)}{p_4(x)},
\]
with $p_i(x)$ polynomials of degree $m$, and then $r_{3m,2m}(x)$ can be computed with $m-1$ products and two inverses.
If $m$ is odd, then $p_2(x)$ and $p_4(x)$ are odd polynomials with complex coefficients.
To avoid this, we consider even $m$ such that $p_2(x)$ and $p_4(x)$ contain pairs of complex conjugate roots, which can be distributed in different ways (and can differ at round off accuracy).
In addition, we have that
\[
p_0(0) + \frac{p_1(0)}{p_2(0)}+ \frac{p_3(0)}{p_4(0)}=1,
\]
and we will take the choice such that $p_0(0)=0$ and $\frac{p_1(0)}{p_2(0)}= \frac{p_3(0)}{p_4(0)}=\frac12$.

Finally, we have also considered the case $r_{4m,3m}(x)$ with $m$ even
\[
r_{4m,3m}(x)=p_0(x) + \frac{p_1(x)}{p_2(x)}+ \frac{p_3(x)}{p_4(x)}+ \frac{p_5(x)}{p_6(x)},
\]
with $p_i(x)$ polynomials of degree $m$, and then $r_{4m,3m}(x)$ can be computed with $m-1$ products and three inverses.
We have analyzed the case $m=4$, i.e. $r_{16,12}(x)$, has a higher order and more accurate results than $r_{13,13}(x)$ while being slightly cheaper.
There is much freedom in the distribution of the roots for $p_{2k}(x), \ k=1,2,3$ (15 different combinations).
In addition, we have that $p_0(0) + \frac{p_1(0)}{p_2(0)}+ \frac{p_3(0)}{p_4(0)}+ \frac{p_5(0)}{p_6(0)}=1$.
The main issue with this scheme is that for all combinations we have explored, the scheme has large positive and negative coefficients which introduce relatively large round off errors so, the scheme is still under investigation.

We have also explored other choices, and the only one we found efficient was the scheme $r_{8,5}(x)$ decomposed as
\[
r_{8,5}(x)=p_0(x) + \frac{p_1(x)}{p_2(x)}+ \frac{q_1(x)}{q_2(x)},
\]
with $q_i(x),p_i(x)$ being polynomials of degree 2 and 3, respectively, and then it can be computed with $2$ products and two inverses.

We followed a similar procedure to reduce the cost of certain diagonal Pad\'e approximants.
We decompose $r_{m,m}(x),\ m=4,8$ into two fractions, while for $m=6,12$ we found more economical to split each one into three fractions.
Different distributions of the roots of the denominators between the fractions have been analyzed to reduce round-off errors.

\subsection{Efficient computation of Taylor polynomial approximations}
\label{sec.4.2}

With $k=0,1,2$ products we can trivially evaluate  $t_n(x)$ for $n=1,2,4$, whereas for $n=8,12$ and $18$ the polynomials can exactly be computed with only $k=3,4$ and 5 products, respectively \cite{preprint,bader19ctm}.

In general, with $k$ products it is possible to build polynomials of degree $2^k$, but for $k>3$ it is not possible to compute $t_{2^k}(x)$ since there are not enough independent parameters.
However, with $k=4,5$ products it is possible to build polynomials of degrees 16 and 24 which coincide with the Taylor expansion up to order 15 and 21 which we will denote by $t_{15}^{[16]}(x)$ and $t_{21}^{[24]}(x)$, respectively.
These methods correspond to the schemes $y_{22}$ and $y_{23}$ in  \cite{sastre19btc}.

As an illustration for $k=3$ products, we show the following algorithm to evaluate $t_8$:
\begin{equation}  \label{Algorithm83eA}
\begin{aligned}
	A_2    & = A^2,                                     \\
	A_4    & = A_2(x_1 A+x_2 A_2),                      \\
	A_8    & = (x_3 A_2+A_4)(x_4I+x_5A+x_6 A_2+x_7A_4), \\
	t_8(A) & = y_0 I + y_1 A +y_{2} A_2+ A_8,
\end{aligned}
\end{equation}
where
\[
\begin{array}{llll}
x_1=\displaystyle x_3\frac{ 1 + \sqrt{177}}{88},           &
x_2= \displaystyle \frac{1 + \sqrt{177}}{352}{x_3},        &
x_4= \displaystyle \frac{-271 + 29\sqrt{177}}{315 x_3},        \\
x_5= \displaystyle \frac{11 (-1 + \sqrt{177})}{1260 x_3},  &
x_6=  \displaystyle \frac{11 (-9 + \sqrt{177})}{5040 x_3}, &
x_7=  \displaystyle \frac{89 - \sqrt{177}}{5040 x_3^2},    &   \\
y_0=1,                                                     &
y_1= 1,                                                    &
y_2 = \displaystyle \frac{857 - 58\sqrt{177}}{630},            \\
x_3 = 2/3.
\end{array}
\]
Notice that, $t_7(x)$ requires at least 4 products, so $t_8(x)$ may be considered a singular polynomial.

\subsection{Selected methods sorted by computational cost}
\label{sec.list}
We collect now the list of the best methods \( w_{\alpha} \) we have found, sorted by their computational cost \( k_{\alpha} \).
The cost to evaluate the inverse is taken as $\frac43\equiv 1\frac13$ products. We give the structure of the schemes with appropriate values for the coefficients $\alpha_i,\beta_i,\ldots, \ i=1,2,\ldots,$ which take different values in each case.
\begin{description}

\item[$k=1$ product] $t_2(x)=1+x+\frac12x^2$.

\item[$k=1\frac13$ product] $r_{2,1}(x)=p_0(x)+\frac{p_1(x)}{p_2(x)}$, $ p_i(x)=\alpha_i+\beta_ix$ with $p_0(0)=0, \alpha_2=1$.

\item[$k=2$ products] $t_4(x)=1+x+x^2(\frac1{2!}+\frac1{3!}x+\frac1{4!}x^2)$.

\item[$k=2\frac13$ products]  $r_{4,2}(x)=p_0(x)+\frac{p_1(x)}{p_2(x)}$, $ p_i(x)=\alpha_i+\beta_ix+\gamma_ix^2$, $\alpha_0=0, \alpha_2=1$.

\item[$k=3$ products] $t_8(x)$.

\item[$k=3\frac13$ products]   $r_{6,3}(x)=p_0(x)+\frac{p_1(x)}{p_2(x)}, \quad p_i(x)=\alpha_i+\beta_ix+\gamma_ix^2+\delta_ix^3$ with $\alpha_0=0, \alpha_2=1$.

\item[$k=3\frac23$ products]   $r_{6,4}(x)=p_0(x)+\frac{p_1(x)}{p_2(x)}+\frac{p_3(x)}{p_4(x)}, \quad p_i(x)=\alpha_i+\beta_ix+\gamma_ix^2$ with $\alpha_0=0$ and $\frac{p_1(0)}{p_2(0)}=\frac{p_3(0)}{p_4(0)}=\frac12, \ \alpha_2=\alpha_4=1$.

\item[$k=4$ products] $t_{12}(x)$ and  $t_{15}^{[16]}(x)$.

\item[$k=4\frac13$ products]   $r_{8,4}(x)=p_0(x)+\frac{p_1(x)}{p_2(x)}, \quad p_i(x)=\alpha_i+\beta_ix+\gamma_ix^2+\delta_ix^3+\sigma_ix^4$ with $\alpha_0=0, \alpha_2=1$.

\item[$k=4\frac23$ products]   $r_{8,5}(x)=p_0(x)+\frac{p_1(x)}{p_2(x)}+\frac{q_1(x)}{q_2(x)}, \quad p_i(x)=\alpha_i+\beta_ix+\gamma_ix^2+\delta_ix^3$, $q_i(x)=\tilde \alpha_i+\tilde \beta_ix+\tilde \gamma_ix^2+\tilde \delta_ix^3$ with $\alpha_0=0$ and $\frac{p_1(0)}{p_2(0)}=\frac{q_1(0)}{q_2(0)}=\frac12$.

\item[$k=5$ products] $t_{18}(x)$ and $t_{21}^{[24]}(x)$.

\item[$k=5\frac13$ products]  $r_{10,5}(x)=p_0(x)+\frac{p_1(x)}{p_2(x)}, \quad p_i(x)=\alpha_i+\beta_ix+\gamma_ix^2+\delta_ix^3+\sigma_ix^4+\mu_ix^5, \ \alpha_0=0, \ \alpha_2=1$.
This method is not used because for all the considered tolerance values its \( \theta \)s  are less than the corresponding ones of \( (r_{8,4}(x))^{2} \).

\item[$k=5\frac23$ products]  $r_{12,8}(x)=p_0(x)+\frac{p_1(x)}{p_2(x)}+\frac{p_3(x)}{p_4(x)}, \quad p_i(x)=\alpha_i+\beta_ix+\gamma_ix^2+\delta_ix^3+\sigma_ix^4$ with $\alpha_0=0$ and $\frac{p_1(0)}{p_2(0)}=\frac{p_3(0)}{p_4(0)}=\frac12, \  \alpha_2=\alpha_4=1$.
Unfortunately, in general it suffers from round-off errors, although smaller than those of $r_{16,12}(x)$.
However for lower tolerance values it can compete with \( r_{8,4}(x) \) and \( r_{8,5}(x) \).

\item[$k=7$ products] $r_{16,12}(x)=p_0(x)+\frac{p_1(x)}{p_2(x)}+\frac{p_3(x)}{p_4(x)}+\frac{p_5(x)}{p_6(x)}, \quad p_i(x)=\alpha_i+\beta_ix+\gamma_ix^2+\delta_ix^3+\sigma_ix^4$. This scheme would be the method of choice when very high accuracy is desired and the scaling has to be used.
Unfortunately, for all choices in the free parameters of the method that we have tried, it suffers severe round-off errors and this scheme is under investigation at this moment and not used in the algorithm.

\item[$k=7\frac13$ products]  $r_{13,13}(x)$.
\end{description}

\paragraph{Comparison of methods for different values of the tolerance.}
In \autoref{tab.thetaAll} and \autoref{tab.thetaDiag} we collect the tolerance values \( tol \) and the corresponding $\theta_{\alpha}(tol)$ for each of the selected methods $w_{\alpha}$ and some additional methods for comparison.
For simplicity in the presentation, the tables show only commonly used values for the tolerance, while for the experimental implementation we use extended tables for \( tol=10^{-k},\ k=0,1,\ldots,16 \).
Methods suboptimal in terms of computational cost for a given tolerance are grayed out and are not used in the final algorithm.
When refining the lookup tables (LUT), we chose to exclude methods with costs lower than that of a squared cheaper method, along with schemes showing large round-off errors.
For example, \( r_{10,5}(x) \) has smaller \( \theta \)s than \( r_{8,4}(x) \) with one additional scaling.

It should be noted that depending on the application, the method search can be implemented in various ways:
\begin{itemize}
	\item If a problem calls for avoiding inverses, the LUT can contain only Tailor methods.
	\item If one expects to work only with high tolerance, higher-order methods like \( r_{13,13}(x) \) may be disabled, thus giving preference to the cheap lower-order methods with additional scaling--squaring.
	\item On the other hand, if one aims to avoid squaring due to low tolerance requirements, squared methods can be additionally penalized, for example, by artificially increasing their cost or by multiplying their \( \theta \) by a factor less than 1.
	\item Additional criteria applied, for instance, the ratio of \( \tfrac{\theta}{k_{\alpha}} \).
\end{itemize}

\paragraph{Algorithm for adaptive scheme selection}
	Taking into account the considerations above, our experimental implementation uses the following algorithm to select a method when provided a matrix, a tolerance value, and a LUT (ordered by tolerance and method cost in ascending order):
\begin{enumerate}
	\item Calculate \( m = \lfloor \log_{10}(tol) \rfloor \), where $\lfloor \cdot\rfloor$ denotes the rounding-down operation.
	\item Select the column with \( \theta_{\alpha} \) corresponding to the condition \( 10^{m} < tol \).
	\item For each \( \theta_{\alpha} \) calculate the corresponding scaling \( s_{\alpha} = \max \left\lbrace 0, \left\lceil \log_{2} \tfrac{\opnorm{A}}{\theta_{\alpha}(10^{m})} \right\rceil \right\rbrace \), where $\lceil \cdot\rceil$ denotes the rounding-up operation.
	\item For each method \( w_{\alpha} \) with the cost \( k_{\alpha} \) (provided in the list above) calculate the total cost \( k_{\alpha} + 1.1 \cdot s_{\alpha} \).
	To break ties, methods that need scaling are penalized by multiplying \( s_{\alpha} \) by the factor of \( 1.1 \).
	\item Select the method with the lowest total cost.
\end{enumerate}

\begin{table}
\caption{\label{tab.thetaAll} Comparison of backward errors of superdiagonal methods.
The asterisk \( (*) \) stands for methods with no cost estimate. Grayed-out rows and columns are provided for completeness.
}
\begin{tblr}{
		rowsep=0.4ex,
		colspec={ll SS SS SS S},
		cell{10,11,16}{1-9} = {fg=lightgray}, %
		cell{18,20}{1-9} = {fg=lightgray}, %
		cell{2-21}{3,5,8} = {fg=lightgray},
	}
	\toprule
	\SetCell[r=2]{c} Cost & \SetCell[r=2]{c} Method & \SetCell[c=7]{c} {{{Tolerance}}} &                 &                 &                 &                  &                 &                  \\ \cmidrule{3-9}
	&                         & {{{$2^{-11}$}}}                  & {{{$10^{-4}$}}} & {{{$2^{-24}$}}} & {{{$10^{-8}$}}} & {{{$10^{-12}$}}} & {{{$2^{-53}$}}} & {{{$10^{-16}$}}} \\ \midrule
	$1$                   & $t_{2}$                 & 5.3053e-2                        & 2.4272e-2       & 5.9789e-4       & 2.4493e-4       & 2.4495e-6        & 2.5810e-8       & 2.4495e-8        \\
	$1\tfrac{1}{3}$       & $r_{2,1}$               & 3.1768e-1                        & 1.8970e-1       & 1.6227e-2       & 8.9557e-3       & 4.1600e-4        & 1.9995e-5       & 1.9310e-5        \\
	$2$                   & $t_{4}$                 & 4.4792e-1                        & 3.1019e-1       & 5.1166e-2       & 3.2872e-2       & 3.3075e-3        & 3.3972e-4       & 3.3095e-4        \\
	$2\tfrac{1}{3}$       & $r_{4,2}$               & 1.6583                           & 1.3026          & 3.9826e-1       & 2.9734e-1       & 6.4820e-2        & 1.4246e-2       & 1.4000e-2        \\
	$3$                   & $t_{8}$                 & 1.5945                           & 1.3454          & 5.8005e-1       & 4.6986e-1       & 1.5397e-1        & 4.9912e-2       & 4.9268e-2        \\
	$3\tfrac{1}{3}$       & $r_{6,3}$               & 3.2781                           & 2.8106          & 1.3146          & 1.0878          & 4.0114e-1        & 1.4715e-1       & 1.4546e-1        \\
	$3\tfrac{2}{3}$       & $r_{6,4}$               & 4.1026                           & 3.5656          & 1.7888          & 1.5071          & 6.1248e-1        & 2.4822e-1       & 2.4565e-1        \\
	$4$                   & $t_{12}$                & 2.7916                           & 2.5021          & 1.4617          & 1.2778          & 6.2401e-1        & 2.9962e-1       & 2.9708e-1        \\
	$*$                   & $t_{15}$                & 3.6842                           & 3.3793          & 2.2170          & 1.9960          & 1.1400           & 6.4108e-1       & 6.3680e-1        \\
	$4$                   & $t_{15}^{[16]}$         & 3.9120                           & 3.5856          & 2.3462          & 2.1113          & 1.2039           & 4.9236e-1       & 4.6327e-1        \\
	$4\tfrac{1}{3}$       & $r_{8,4}$               & 4.9543                           & 4.4284          & 2.5478          & 2.2191          & 1.0668           & 5.0739e-1       & 5.0305e-1        \\
	$4\tfrac{2}{3}$       & $r_{8,5}$               & 5.8331                           & 5.2529          & 3.1401          & 2.7621          & 1.4012           & 7.0491e-1       & 6.9934e-1        \\
	$5$                   & $t_{18}$                & 4.5703                           & 4.2556          & 3.0101          & 2.7620          & 1.7473           & 1.0909          & 1.0849           \\
	$*$                   & $t_{21}$                & 5.4505                           & 5.1293          & 3.8239          & 3.5557          & 2.4160           & 1.6237          & 1.6162           \\
	$5$                   & $t_{21}^{[24]}$         & 5.6233                           & 5.2926          & 3.9496          & 3.6737          & 2.4998           & 4.5420e-1       & 4.2091e-1        \\
	$5\tfrac{1}{3}$       & $r_{10,5}$              & 6.6426                           & 6.0821          & 3.9474          & 3.5435          & 1.9959           & 1.1108          & 1.1033           \\
	$5\tfrac{2}{3}$       & $r_{12,8}$              & 1.0200e1                         & 9.5441          & 6.9059          & 6.3724          & 4.1589           & 2.6901          & 2.6765           \\
	$7$                   & $r_{16,12}$             & 1.5542e1                         & 1.4825e1        & 1.1799e1        & 1.1152e1        & 8.2701           & 6.0934          & 6.0718           \\
	$7\tfrac{1}{3}$       & $r_{13,13}$             & 1.5331e1                         & 1.4542e1        & 1.1249e1        & 1.0557e1        & 7.5495           & 5.3719          & 5.3508           \\ \bottomrule
\end{tblr}
\end{table}

\section{The matrix exponential in Lie groups}
Given $J\in GL(N)$ with $GL(N)$ being the group of all $N\times N$ non-singular matrices, one can define the $J$-orthogonal group, or quadratic Lie group, as the set
\[
O_J(N)=\{ X\in GL(N): X^{\dagger}JX=J\},
\]
($X^{\dagger}$ stands for conjugate transpose)
and the Lie algebra of the quadratic Lie group is given by
\[
o_J(N)=\{ A\in gl(N): A^{\dagger}J+JA=0_N \},
\]
where $gl(N)$ is the algebra of all $N\times N$  matrices
with $0_N$ the null marix of dimension $N$.

It is well known that if $A\in o_J(N)$ then $\e^A\in O_J(N)$. If $J$ is the identity matrix this is satisfied when $A$ is a skew-Hermitian (skew-symmetric) matrix and the exponential matrix is a unitary (orthogonal) matrix while if $J$ is the symplectic matrix then $\e^A$ belongs to the symplectic group, etc.

Most Lie group integrators are of relatively low order and usually provide moderate accuracy (for example, \eqref{eq.Crouch-Grossman1}), but for many problems it is essential to preserve the Lie group structure when calculating the exponential \cite{celledoni00ate,celledoni01mft,iserles00lgm}.
In such case, among the schemes considered in this work, one should use diagonal Pad\'e approximants which are the only ones that exactly preserve such structure.
Notice that, since $A^{\dagger}J=-JA$ then
\[
\left(A^{\dagger}\right)^kJ=(-1)^kJA^k,
\]
so, if $k$ is odd then $A^k$ belongs to the algebra and its exponential belongs to the Lie group. A similar property occurs for any odd function $f(A)$, i.e. if $f(-A)=-f(A)$ where $f(x)$ is an analytic function at the origin then $f(A)$ belongs to the Lie algebra and $\e^{f(A)}$ belongs to the Lie group.
Diagonal Pad\'e approximants are the only Pad\'e methods which are symmetric, so
\[
r_{m,m}(-x)=(r_{m,m}(x))^{-1},
\]
or, equivalently,
\[
r_{m,m}(x)=\e^{f(x)} \qquad \mbox{with} \qquad f(-x)=-f(x),
\]
and then they preserve the Lie group structure.

Our algorithm has an option in order to approximate the matrix exponential using diagonal Pad\'e approximants.
In this case we will use the schemes already incorporated in $\textsc{MATLAB}$, but we will add some few more schemes following the fractional decomposition previously mentioned. They are collected in \autoref{tab.thetaDiag}.
Note that methods with \( m=10,\ 11 \) would not be used because \( r_{12,12}(x) \) has higher order with the same computational cost.
Unfortunately, the decomposed  \( r_{12,12}(x) \) suffers from the same round-off errors as \( r_{12,8}(x) \) and \( r_{16,12}(x) \), so it is not used in the algorithm.

\begin{table}
\caption{\label{tab.thetaDiag} Comparison of backward errors of diagonal Padé approximants.
Grayed-out rows and columns are provided for completeness.
}
\begin{tblr}{
		rowsep=0.4ex,
		colspec={llSSSSSSS},
		cell{11-12}{1-9} = {fg=lightgray},%
		cell{15-20}{1-9} = {fg=lightgray},%
		cell{2-21}{3,5,8} = {fg=lightgray},
	}
	\toprule
	\SetCell[r=2]{c} Cost & \SetCell[r=2]{c} Method & \SetCell[c=7]{c} {{{Tolerance}}} &                 &                 &                 &                  &                 &                  \\ \cmidrule{3-9}
	&                         & {{{$2^{-11}$}}}                  & {{{$10^{-4}$}}} & {{{$2^{-24}$}}} & {{{$10^{-8}$}}} & {{{$10^{-12}$}}} & {{{$2^{-53}$}}} & {{{$10^{-16}$}}} \\ \midrule
	$2\frac13$            & $r_{2,2}$               & 7.6339e-1                        & 5.1596e-1       & 8.0930e-2       & 5.1798e-2       & 5.1800e-3        & 5.3172e-4       & 5.1800e-4        \\
	$3\frac13$            & $r_{3,3}$               & 1.8718                           & 1.4500          & 4.2587e-1       & 3.1644e-1       & 6.8218e-2        & 1.4956e-2       & 1.4697e-2        \\
	$3\frac23$            & $r_{4,4}$               & 3.1358                           & 2.6004          & 1.0490          & 8.4041e-1       & 2.6638e-1        & 8.5364e-2       & 8.4255e-2        \\
	$4\frac13$            & $r_{5,5}$               & 4.4590                           & 3.8495          & 1.8802          & 1.5766          & 6.3074e-1        & 2.5394e-1       & 2.5130e-1        \\
	$5$                   & $r_{6,6}$               & 5.8066                           & 5.1466          & 2.8543          & 2.4680          & 1.1545           & 5.4147e-1       & 5.3677e-1        \\
	$5\frac13$            & $r_{7,7}$               & 7.1643                           & 6.4685          & 3.9257          & 3.4697          & 1.8161           & 9.5042e-1       & 9.4336e-1        \\
	$5\frac23$            & $r_{8,8}$               & 8.5260                           & 7.8037          & 5.0640          & 4.5498          & 2.5917           & 1.4732          & 1.4636           \\
	$6\frac13$            & $r_{9,9}$               & 9.8887                           & 9.1462          & 6.2492          & 5.6866          & 3.4599           & 2.0978          & 2.0858           \\
	7                     & $r_{10,10}$             & 1.1251e1                         & 1.0493e1        & 7.4679          & 6.8648          & 4.4027           & 2.8116          & 2.7971           \\
	7                     & $r_{11,11}$             & 1.2613e1                         & 1.1842e1        & 8.7113          & 8.0739          & 5.4058           & 3.6023          & 3.5855           \\
	$7$                   & $r_{12,12}$             & 1.3973e1                         & 1.3191e1        & 9.9731          & 9.3064          & 6.4578           & 4.4589          & 4.4399           \\
	$7\frac13$            & $r_{13,13}$             & 1.5331e1                         & 1.4542e1        & 1.1249e1        & 1.0557e1        & 7.5495           & 5.3719          & 5.3508           \\
	$8\frac13$            & $r_{14,14}$             & 1.6689e1                         & 1.5892e1        & 1.2535e1        & 1.1822          & 8.6739           & 6.3331          & 6.3101           \\
	$8\frac13$            & $r_{15,15}$             & 1.8045e1                         & 1.7242e1        & 1.3830e1        & 1.3097e1        & 9.8254           & 7.3357          & 7.3109           \\
	$8\frac13$            & $r_{16,16}$             & 1.9399e1                         & 1.8591e1        & 1.5131e1        & 1.4382e1        & 1.0999e1         & 8.3737          & 8.3473           \\
	9                     & $r_{17,17}$             & 2.0752e1                         & 1.9940e1        & 1.6438e1        & 1.5674e1        & 1.2192e1         & 9.4424          & 9.4144           \\
	9                     & $r_{18,18}$             & 2.2105e1                         & 2.1288e1        & 1.7749e1        & 1.6972e1        & 1.3401e1         & 1.0537e1        & 1.0508e1         \\
	\bottomrule
\end{tblr}
\end{table}

\section{Numerical performance}
\label{sec.5}
In the numerical experiments we show the behavior of the proposed algorithm and the individual schemes that comprise it.
In the examples, the algorithm calculates \( \opnorm{A} \) and uses it along with the prescribed tolerance \( tol \) to select the most suitable scheme from \autoref{tab.thetaAll} or \autoref{tab.thetaDiag} as described in \autoref{sec.list}.

\paragraph{Numerical example 1} Consider a random diagonally dominant matrix \( A = D + R, \ A \in \mathbb{R}^{101 \times 101} \), such that \( D = diag(-50,-49,\ldots,50) \) and the elements of \( R \) are sampled uniformly from \( \left[ -1;\,1 \right] \), and the matrix is normalized (i.e. $A:=A/\opnorm{A}$).

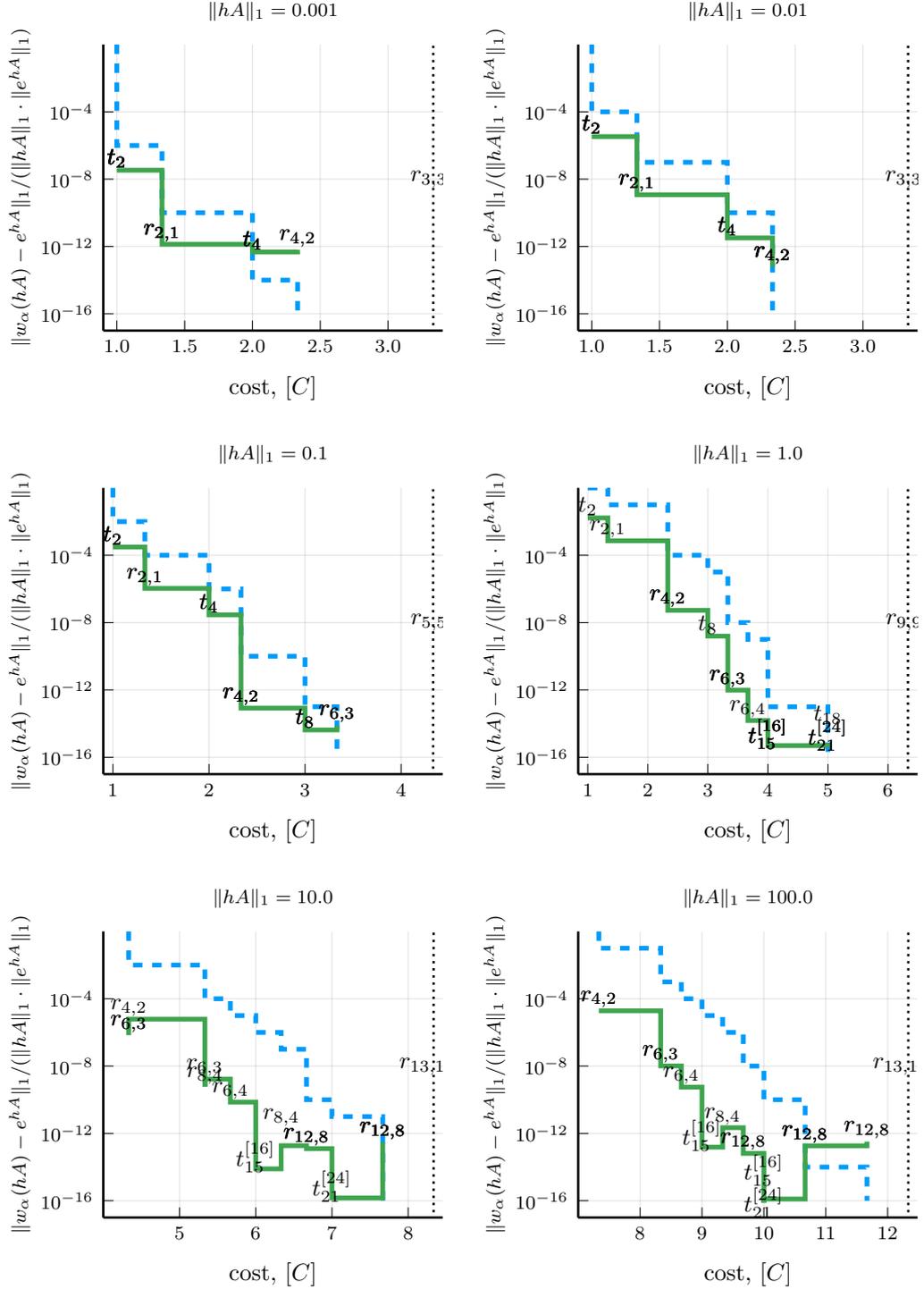
\begin{figure}[p]
\begin{subfigure}[b]{0.48\textwidth}
	\raggedleft
	\begin{tikzpicture}[/tikz/background rectangle/.style={fill={rgb,1:red,1.0;green,1.0;blue,1.0}, fill opacity={1.0}, draw opacity={1.0}}, show background rectangle]
\begin{axis}[point meta max={nan}, point meta min={nan}, legend cell align={left}, legend columns={1}, title={$\lVert hA \rVert_{1}=0.001$}, title style={at={{(0.5,1)}}, anchor={south}, color={rgb,1:red,0.0;green,0.0;blue,0.0}, draw opacity={1.0}, rotate={0.0}, align={center}}, legend style={color={rgb,1:red,0.0;green,0.0;blue,0.0}, draw opacity={1.0}, line width={1}, solid, fill={rgb,1:red,1.0;green,1.0;blue,1.0}, fill opacity={1.0}, text opacity={1.0}, text={rgb,1:red,0.0;green,0.0;blue,0.0}, cells={anchor={center}}, at={(1.02, 1)}, anchor={north west}}, axis background/.style={fill={rgb,1:red,1.0;green,1.0;blue,1.0}, opacity={1.0}}, anchor={north west}, scaled x ticks={false}, xlabel={cost,\ [$C$]}, x tick style={color={rgb,1:red,0.0;green,0.0;blue,0.0}, opacity={1.0}}, x tick label style={color={rgb,1:red,0.0;green,0.0;blue,0.0}, opacity={1.0}, rotate={0}}, xmajorgrids={true}, xmin={0.9}, xmax={3.4033333333333338}, xticklabels={{$1.0$,$1.5$,$2.0$,$2.5$,$3.0$}}, xtick={{1.0,1.5,2.0,2.5,3.0}}, xtick align={inside}, xticklabel style={color={rgb,1:red,0.0;green,0.0;blue,0.0}, draw opacity={1.0}, rotate={0.0}}, x grid style={color={rgb,1:red,0.0;green,0.0;blue,0.0}, draw opacity={0.1}, line width={0.5}, solid}, axis x line*={left}, x axis line style={color={rgb,1:red,0.0;green,0.0;blue,0.0}, draw opacity={1.0}, line width={1}, solid}, scaled y ticks={false}, ylabel={$\lVert w_{\alpha}(h A)-e^{h A} \rVert_{1} / (\lVert hA \rVert_{1} \cdot \lVert e^{hA} \rVert_{1})$}, y tick style={color={rgb,1:red,0.0;green,0.0;blue,0.0}, opacity={1.0}}, y tick label style={color={rgb,1:red,0.0;green,0.0;blue,0.0}, opacity={1.0}, rotate={0}}, ymode={log}, log basis y={10}, ymajorgrids={true}, ymin={1.0e-17}, ymax={1.0}, yticklabels={{$10^{-16}$,$10^{-8}$,$10^{-12}$,$10^{-4}$}}, ytick={{1.0e-16,1.0e-8,1.0e-12,0.0001}}, ytick align={inside}, yticklabel style={color={rgb,1:red,0.0;green,0.0;blue,0.0}, draw opacity={1.0}, rotate={0.0}}, y grid style={color={rgb,1:red,0.0;green,0.0;blue,0.0}, draw opacity={0.1}, line width={0.5}, solid}, axis y line*={left}, y axis line style={color={rgb,1:red,0.0;green,0.0;blue,0.0}, draw opacity={1.0}, line width={1}, solid}, colorbar={false}]
    \addplot[color={rgb,1:red,0.0;green,0.6056;blue,0.9787}, name path={6679839a-797a-472c-913a-8807be60a45b}, const plot, draw opacity={1.0}, line width={2}, dashed]
        table[row sep={\\}]
        {
            \\
            1.0  1.0  \\
            1.0  0.1  \\
            1.0  0.010000000000000002  \\
            1.0  0.001  \\
            1.0  0.0001  \\
            1.0  1.0e-5  \\
            1.0  1.0e-6  \\
            1.3333333333333317  1.0e-7  \\
            1.3333333333333317  1.0e-8  \\
            1.3333333333333317  1.0e-9  \\
            1.3333333333333317  1.0e-10  \\
            2.0  1.0e-11  \\
            2.0  1.0e-12  \\
            2.0  1.0e-13  \\
            2.0  1.0e-14  \\
            2.3333333333333366  1.0e-15  \\
            2.3333333333333366  1.0e-16  \\
        }
        ;
    \addplot[color={rgb,1:red,0.0;green,0.0;blue,0.0}, name path={9ac4a69e-9b20-4b35-827f-01e6676def7e}, draw opacity={1.0}, line width={1}, dotted]
        table[row sep={\\}]
        {
            \\
            3.3333333333333335  1.0e-34  \\
            3.3333333333333335  1.0e17  \\
        }
        ;
    \addplot[color={rgb,1:red,0.2422;green,0.6433;blue,0.3044}, name path={a9f74102-545c-4c29-9fad-102ca806797a}, const plot, draw opacity={1.0}, line width={2}, solid]
        table[row sep={\\}]
        {
            \\
            1.0  3.3784368133331593e-8  \\
            1.0  3.3784368133331593e-8  \\
            1.0  3.3784368133331593e-8  \\
            1.0  3.3784368133331593e-8  \\
            1.0  3.3784368133331593e-8  \\
            1.0  3.3784368133331593e-8  \\
            1.0  3.3784368133331593e-8  \\
            1.3333333333333317  1.3450077453861433e-12  \\
            1.3333333333333317  1.3450077453861433e-12  \\
            1.3333333333333317  1.3450077453861433e-12  \\
            1.3333333333333317  1.3450077453861433e-12  \\
            2.0  4.671564322735404e-13  \\
            2.0  4.671564322735404e-13  \\
            2.0  4.671564322735404e-13  \\
            2.0  4.671564322735404e-13  \\
            2.3333333333333366  6.257676594015149e-13  \\
            2.3333333333333366  6.257676594015149e-13  \\
        }
        ;
    \node[]  at (axis cs:3.3000000000000003,1.0e-8) {$r_{3,3}$};
    \node[]  at (axis cs:1.0,1.6892184066665798e-7) {$t_{2}$};
    \node[]  at (axis cs:1.0,1.6892184066665798e-7) {$t_{2}$};
    \node[]  at (axis cs:1.0,1.6892184066665798e-7) {$t_{2}$};
    \node[]  at (axis cs:1.0,1.6892184066665798e-7) {$t_{2}$};
    \node[]  at (axis cs:1.0,1.6892184066665798e-7) {$t_{2}$};
    \node[]  at (axis cs:1.0,1.6892184066665798e-7) {$t_{2}$};
    \node[]  at (axis cs:1.0,1.6892184066665798e-7) {$t_{2}$};
    \node[]  at (axis cs:1.3333333333333317,6.7250387269307166e-12) {$r_{2,1}$};
    \node[]  at (axis cs:1.3333333333333317,6.7250387269307166e-12) {$r_{2,1}$};
    \node[]  at (axis cs:1.3333333333333317,6.7250387269307166e-12) {$r_{2,1}$};
    \node[]  at (axis cs:1.3333333333333317,6.7250387269307166e-12) {$r_{2,1}$};
    \node[]  at (axis cs:2.0,2.3357821613677024e-12) {$t_{4}$};
    \node[]  at (axis cs:2.0,2.3357821613677024e-12) {$t_{4}$};
    \node[]  at (axis cs:2.0,2.3357821613677024e-12) {$t_{4}$};
    \node[]  at (axis cs:2.0,2.3357821613677024e-12) {$t_{4}$};
    \node[]  at (axis cs:2.3333333333333366,3.128838297007575e-12) {$r_{4,2}$};
    \node[]  at (axis cs:2.3333333333333366,3.128838297007575e-12) {$r_{4,2}$};
\end{axis}
\end{tikzpicture}%
\end{subfigure}
~
\begin{subfigure}[b]{0.48\textwidth}
	\raggedright
	\begin{tikzpicture}[/tikz/background rectangle/.style={fill={rgb,1:red,1.0;green,1.0;blue,1.0}, fill opacity={1.0}, draw opacity={1.0}}, show background rectangle]
\begin{axis}[point meta max={nan}, point meta min={nan}, legend cell align={left}, legend columns={1}, title={$\lVert hA \rVert_{1}=0.01$}, title style={at={{(0.5,1)}}, anchor={south}, color={rgb,1:red,0.0;green,0.0;blue,0.0}, draw opacity={1.0}, rotate={0.0}, align={center}}, legend style={color={rgb,1:red,0.0;green,0.0;blue,0.0}, draw opacity={1.0}, line width={1}, solid, fill={rgb,1:red,1.0;green,1.0;blue,1.0}, fill opacity={1.0}, text opacity={1.0}, text={rgb,1:red,0.0;green,0.0;blue,0.0}, cells={anchor={center}}, at={(1.02, 1)}, anchor={north west}}, axis background/.style={fill={rgb,1:red,1.0;green,1.0;blue,1.0}, opacity={1.0}}, anchor={north west}, scaled x ticks={false}, xlabel={cost,\ [$C$]}, x tick style={color={rgb,1:red,0.0;green,0.0;blue,0.0}, opacity={1.0}}, x tick label style={color={rgb,1:red,0.0;green,0.0;blue,0.0}, opacity={1.0}, rotate={0}}, xmajorgrids={true}, xmin={0.9}, xmax={3.4033333333333338}, xticklabels={{$1.0$,$1.5$,$2.0$,$2.5$,$3.0$}}, xtick={{1.0,1.5,2.0,2.5,3.0}}, xtick align={inside}, xticklabel style={color={rgb,1:red,0.0;green,0.0;blue,0.0}, draw opacity={1.0}, rotate={0.0}}, x grid style={color={rgb,1:red,0.0;green,0.0;blue,0.0}, draw opacity={0.1}, line width={0.5}, solid}, axis x line*={left}, x axis line style={color={rgb,1:red,0.0;green,0.0;blue,0.0}, draw opacity={1.0}, line width={1}, solid}, scaled y ticks={false}, ylabel={$\lVert w_{\alpha}(h A)-e^{h A} \rVert_{1} / (\lVert hA \rVert_{1} \cdot \lVert e^{hA} \rVert_{1})$}, y tick style={color={rgb,1:red,0.0;green,0.0;blue,0.0}, opacity={1.0}}, y tick label style={color={rgb,1:red,0.0;green,0.0;blue,0.0}, opacity={1.0}, rotate={0}}, ymode={log}, log basis y={10}, ymajorgrids={true}, ymin={1.0e-17}, ymax={1.0}, yticklabels={{$10^{-16}$,$10^{-8}$,$10^{-12}$,$10^{-4}$}}, ytick={{1.0e-16,1.0e-8,1.0e-12,0.0001}}, ytick align={inside}, yticklabel style={color={rgb,1:red,0.0;green,0.0;blue,0.0}, draw opacity={1.0}, rotate={0.0}}, y grid style={color={rgb,1:red,0.0;green,0.0;blue,0.0}, draw opacity={0.1}, line width={0.5}, solid}, axis y line*={left}, y axis line style={color={rgb,1:red,0.0;green,0.0;blue,0.0}, draw opacity={1.0}, line width={1}, solid}, colorbar={false}]
    \addplot[color={rgb,1:red,0.0;green,0.6056;blue,0.9787}, name path={180bf8cf-dc1c-43d3-a843-b13cde8a70d5}, const plot, draw opacity={1.0}, line width={2}, dashed]
        table[row sep={\\}]
        {
            \\
            1.0  1.0  \\
            1.0  0.1  \\
            1.0  0.010000000000000002  \\
            1.0  0.001  \\
            1.0  0.0001  \\
            1.3333333333333317  1.0e-5  \\
            1.3333333333333317  1.0e-6  \\
            1.3333333333333317  1.0e-7  \\
            2.0  1.0e-8  \\
            2.0  1.0e-9  \\
            2.0  1.0e-10  \\
            2.3333333333333366  1.0e-11  \\
            2.3333333333333366  1.0e-12  \\
            2.3333333333333366  1.0e-13  \\
            2.3333333333333366  1.0e-14  \\
            2.3333333333333366  1.0e-15  \\
            2.3333333333333366  1.0e-16  \\
        }
        ;
    \addplot[color={rgb,1:red,0.0;green,0.0;blue,0.0}, name path={70c18779-0d31-4ea6-87ee-efed9067f5a3}, draw opacity={1.0}, line width={1}, dotted]
        table[row sep={\\}]
        {
            \\
            3.3333333333333335  1.0e-34  \\
            3.3333333333333335  1.0e17  \\
        }
        ;
    \addplot[color={rgb,1:red,0.2422;green,0.6433;blue,0.3044}, name path={b59eedbf-b887-4046-bf60-6c6b33a485e6}, const plot, draw opacity={1.0}, line width={2}, solid]
        table[row sep={\\}]
        {
            \\
            1.0  3.382640252824104e-6  \\
            1.0  3.382640252824104e-6  \\
            1.0  3.382640252824104e-6  \\
            1.0  3.382640252824104e-6  \\
            1.0  3.382640252824104e-6  \\
            1.3333333333333317  1.1814328012552832e-9  \\
            1.3333333333333317  1.1814328012552832e-9  \\
            1.3333333333333317  1.1814328012552832e-9  \\
            2.0  3.203194874640667e-12  \\
            2.0  3.203194874640667e-12  \\
            2.0  3.203194874640667e-12  \\
            2.3333333333333366  6.310014072341919e-14  \\
            2.3333333333333366  6.310014072341919e-14  \\
            2.3333333333333366  6.310014072341919e-14  \\
            2.3333333333333366  6.310014072341919e-14  \\
            2.3333333333333366  6.310014072341919e-14  \\
            2.3333333333333366  6.310014072341919e-14  \\
        }
        ;
    \node[]  at (axis cs:3.3000000000000003,1.0e-8) {$r_{3,3}$};
    \node[]  at (axis cs:1.0,1.691320126412052e-5) {$t_{2}$};
    \node[]  at (axis cs:1.0,1.691320126412052e-5) {$t_{2}$};
    \node[]  at (axis cs:1.0,1.691320126412052e-5) {$t_{2}$};
    \node[]  at (axis cs:1.0,1.691320126412052e-5) {$t_{2}$};
    \node[]  at (axis cs:1.0,1.691320126412052e-5) {$t_{2}$};
    \node[]  at (axis cs:1.3333333333333317,5.9071640062764154e-9) {$r_{2,1}$};
    \node[]  at (axis cs:1.3333333333333317,5.9071640062764154e-9) {$r_{2,1}$};
    \node[]  at (axis cs:1.3333333333333317,5.9071640062764154e-9) {$r_{2,1}$};
    \node[]  at (axis cs:2.0,1.6015974373203335e-11) {$t_{4}$};
    \node[]  at (axis cs:2.0,1.6015974373203335e-11) {$t_{4}$};
    \node[]  at (axis cs:2.0,1.6015974373203335e-11) {$t_{4}$};
    \node[]  at (axis cs:2.3333333333333366,3.1550070361709595e-13) {$r_{4,2}$};
    \node[]  at (axis cs:2.3333333333333366,3.1550070361709595e-13) {$r_{4,2}$};
    \node[]  at (axis cs:2.3333333333333366,3.1550070361709595e-13) {$r_{4,2}$};
    \node[]  at (axis cs:2.3333333333333366,3.1550070361709595e-13) {$r_{4,2}$};
    \node[]  at (axis cs:2.3333333333333366,3.1550070361709595e-13) {$r_{4,2}$};
    \node[]  at (axis cs:2.3333333333333366,3.1550070361709595e-13) {$r_{4,2}$};
\end{axis}
\end{tikzpicture}%
\end{subfigure}

\begin{subfigure}[b]{0.48\textwidth}
	\raggedleft
	\begin{tikzpicture}[/tikz/background rectangle/.style={fill={rgb,1:red,1.0;green,1.0;blue,1.0}, fill opacity={1.0}, draw opacity={1.0}}, show background rectangle]
\begin{axis}[point meta max={nan}, point meta min={nan}, legend cell align={left}, legend columns={1}, title={$\lVert hA \rVert_{1}=0.1$}, title style={at={{(0.5,1)}}, anchor={south}, color={rgb,1:red,0.0;green,0.0;blue,0.0}, draw opacity={1.0}, rotate={0.0}, align={center}}, legend style={color={rgb,1:red,0.0;green,0.0;blue,0.0}, draw opacity={1.0}, line width={1}, solid, fill={rgb,1:red,1.0;green,1.0;blue,1.0}, fill opacity={1.0}, text opacity={1.0}, text={rgb,1:red,0.0;green,0.0;blue,0.0}, cells={anchor={center}}, at={(1.02, 1)}, anchor={north west}}, axis background/.style={fill={rgb,1:red,1.0;green,1.0;blue,1.0}, opacity={1.0}}, anchor={north west}, scaled x ticks={false}, xlabel={cost,\ [$C$]}, x tick style={color={rgb,1:red,0.0;green,0.0;blue,0.0}, opacity={1.0}}, x tick label style={color={rgb,1:red,0.0;green,0.0;blue,0.0}, opacity={1.0}, rotate={0}}, xmajorgrids={true}, xmin={0.8999999999999999}, xmax={4.433333333333334}, xticklabels={{$1$,$2$,$3$,$4$}}, xtick={{1.0,2.0,3.0,4.0}}, xtick align={inside}, xticklabel style={color={rgb,1:red,0.0;green,0.0;blue,0.0}, draw opacity={1.0}, rotate={0.0}}, x grid style={color={rgb,1:red,0.0;green,0.0;blue,0.0}, draw opacity={0.1}, line width={0.5}, solid}, axis x line*={left}, x axis line style={color={rgb,1:red,0.0;green,0.0;blue,0.0}, draw opacity={1.0}, line width={1}, solid}, scaled y ticks={false}, ylabel={$\lVert w_{\alpha}(h A)-e^{h A} \rVert_{1} / (\lVert hA \rVert_{1} \cdot \lVert e^{hA} \rVert_{1})$}, y tick style={color={rgb,1:red,0.0;green,0.0;blue,0.0}, opacity={1.0}}, y tick label style={color={rgb,1:red,0.0;green,0.0;blue,0.0}, opacity={1.0}, rotate={0}}, ymode={log}, log basis y={10}, ymajorgrids={true}, ymin={1.0e-17}, ymax={1.0}, yticklabels={{$10^{-16}$,$10^{-8}$,$10^{-12}$,$10^{-4}$}}, ytick={{1.0e-16,1.0e-8,1.0e-12,0.0001}}, ytick align={inside}, yticklabel style={color={rgb,1:red,0.0;green,0.0;blue,0.0}, draw opacity={1.0}, rotate={0.0}}, y grid style={color={rgb,1:red,0.0;green,0.0;blue,0.0}, draw opacity={0.1}, line width={0.5}, solid}, axis y line*={left}, y axis line style={color={rgb,1:red,0.0;green,0.0;blue,0.0}, draw opacity={1.0}, line width={1}, solid}, colorbar={false}]
    \addplot[color={rgb,1:red,0.0;green,0.6056;blue,0.9787}, name path={0219ce9f-1860-4df2-aa63-47e0ac6638e3}, const plot, draw opacity={1.0}, line width={2}, dashed]
        table[row sep={\\}]
        {
            \\
            1.0  1.0  \\
            1.0  0.1  \\
            1.0  0.010000000000000002  \\
            1.3333333333333317  0.001  \\
            1.3333333333333317  0.0001  \\
            2.0  1.0e-5  \\
            2.0  1.0e-6  \\
            2.3333333333333366  1.0e-7  \\
            2.3333333333333366  1.0e-8  \\
            2.3333333333333366  1.0e-9  \\
            2.3333333333333366  1.0e-10  \\
            3.0  1.0e-11  \\
            3.0  1.0e-12  \\
            3.0  1.0e-13  \\
            3.333333333333332  1.0e-14  \\
            3.333333333333332  1.0e-15  \\
            3.333333333333332  1.0e-16  \\
        }
        ;
    \addplot[color={rgb,1:red,0.0;green,0.0;blue,0.0}, name path={00e49165-751f-44a7-8d56-5494af49fc58}, draw opacity={1.0}, line width={1}, dotted]
        table[row sep={\\}]
        {
            \\
            4.333333333333333  1.0e-34  \\
            4.333333333333333  1.0e17  \\
        }
        ;
    \addplot[color={rgb,1:red,0.2422;green,0.6433;blue,0.3044}, name path={52df5560-926b-4ed1-9dcd-941022764275}, const plot, draw opacity={1.0}, line width={2}, solid]
        table[row sep={\\}]
        {
            \\
            1.0  0.000307529056543777  \\
            1.0  0.000307529056543777  \\
            1.0  0.000307529056543777  \\
            1.3333333333333317  1.0752313693126088e-6  \\
            1.3333333333333317  1.0752313693126088e-6  \\
            2.0  2.870367036776465e-8  \\
            2.0  2.870367036776465e-8  \\
            2.3333333333333366  8.351621665666424e-14  \\
            2.3333333333333366  8.351621665666424e-14  \\
            2.3333333333333366  8.351621665666424e-14  \\
            2.3333333333333366  8.351621665666424e-14  \\
            3.0  4.18897824108267e-15  \\
            3.0  4.18897824108267e-15  \\
            3.0  4.18897824108267e-15  \\
            3.333333333333332  7.806133947473576e-15  \\
            3.333333333333332  7.806133947473576e-15  \\
            3.333333333333332  7.806133947473576e-15  \\
        }
        ;
    \node[]  at (axis cs:4.29,1.0e-8) {$r_{5,5}$};
    \node[]  at (axis cs:1.0,0.001537645282718885) {$t_{2}$};
    \node[]  at (axis cs:1.0,0.001537645282718885) {$t_{2}$};
    \node[]  at (axis cs:1.0,0.001537645282718885) {$t_{2}$};
    \node[]  at (axis cs:1.3333333333333317,5.376156846563044e-6) {$r_{2,1}$};
    \node[]  at (axis cs:1.3333333333333317,5.376156846563044e-6) {$r_{2,1}$};
    \node[]  at (axis cs:2.0,1.4351835183882324e-7) {$t_{4}$};
    \node[]  at (axis cs:2.0,1.4351835183882324e-7) {$t_{4}$};
    \node[]  at (axis cs:2.3333333333333366,4.175810832833212e-13) {$r_{4,2}$};
    \node[]  at (axis cs:2.3333333333333366,4.175810832833212e-13) {$r_{4,2}$};
    \node[]  at (axis cs:2.3333333333333366,4.175810832833212e-13) {$r_{4,2}$};
    \node[]  at (axis cs:2.3333333333333366,4.175810832833212e-13) {$r_{4,2}$};
    \node[]  at (axis cs:3.0,2.094489120541335e-14) {$t_{8}$};
    \node[]  at (axis cs:3.0,2.094489120541335e-14) {$t_{8}$};
    \node[]  at (axis cs:3.0,2.094489120541335e-14) {$t_{8}$};
    \node[]  at (axis cs:3.333333333333332,3.9030669737367875e-14) {$r_{6,3}$};
    \node[]  at (axis cs:3.333333333333332,3.9030669737367875e-14) {$r_{6,3}$};
    \node[]  at (axis cs:3.333333333333332,3.9030669737367875e-14) {$r_{6,3}$};
\end{axis}
\end{tikzpicture}%
\end{subfigure}
~
\begin{subfigure}[b]{0.48\textwidth}
	\raggedright
	\begin{tikzpicture}[/tikz/background rectangle/.style={fill={rgb,1:red,1.0;green,1.0;blue,1.0}, fill opacity={1.0}, draw opacity={1.0}}, show background rectangle]
\begin{axis}[point meta max={nan}, point meta min={nan}, legend cell align={left}, legend columns={1}, title={$\lVert hA \rVert_{1}=1.0$}, title style={at={{(0.5,1)}}, anchor={south}, color={rgb,1:red,0.0;green,0.0;blue,0.0}, draw opacity={1.0}, rotate={0.0}, align={center}}, legend style={color={rgb,1:red,0.0;green,0.0;blue,0.0}, draw opacity={1.0}, line width={1}, solid, fill={rgb,1:red,1.0;green,1.0;blue,1.0}, fill opacity={1.0}, text opacity={1.0}, text={rgb,1:red,0.0;green,0.0;blue,0.0}, cells={anchor={center}}, at={(1.02, 1)}, anchor={north west}}, axis background/.style={fill={rgb,1:red,1.0;green,1.0;blue,1.0}, opacity={1.0}}, anchor={north west}, scaled x ticks={false}, xlabel={cost,\ [$C$]}, x tick style={color={rgb,1:red,0.0;green,0.0;blue,0.0}, opacity={1.0}}, x tick label style={color={rgb,1:red,0.0;green,0.0;blue,0.0}, opacity={1.0}, rotate={0}}, xmajorgrids={true}, xmin={0.8399999999999999}, xmax={6.493333333333333}, xticklabels={{$1$,$2$,$3$,$4$,$5$,$6$}}, xtick={{1.0,2.0,3.0,4.0,5.0,6.0}}, xtick align={inside}, xticklabel style={color={rgb,1:red,0.0;green,0.0;blue,0.0}, draw opacity={1.0}, rotate={0.0}}, x grid style={color={rgb,1:red,0.0;green,0.0;blue,0.0}, draw opacity={0.1}, line width={0.5}, solid}, axis x line*={left}, x axis line style={color={rgb,1:red,0.0;green,0.0;blue,0.0}, draw opacity={1.0}, line width={1}, solid}, scaled y ticks={false}, ylabel={$\lVert w_{\alpha}(h A)-e^{h A} \rVert_{1} / (\lVert hA \rVert_{1} \cdot \lVert e^{hA} \rVert_{1})$}, y tick style={color={rgb,1:red,0.0;green,0.0;blue,0.0}, opacity={1.0}}, y tick label style={color={rgb,1:red,0.0;green,0.0;blue,0.0}, opacity={1.0}, rotate={0}}, ymode={log}, log basis y={10}, ymajorgrids={true}, ymin={1.0e-17}, ymax={1.0}, yticklabels={{$10^{-16}$,$10^{-8}$,$10^{-12}$,$10^{-4}$}}, ytick={{1.0e-16,1.0e-8,1.0e-12,0.0001}}, ytick align={inside}, yticklabel style={color={rgb,1:red,0.0;green,0.0;blue,0.0}, draw opacity={1.0}, rotate={0.0}}, y grid style={color={rgb,1:red,0.0;green,0.0;blue,0.0}, draw opacity={0.1}, line width={0.5}, solid}, axis y line*={left}, y axis line style={color={rgb,1:red,0.0;green,0.0;blue,0.0}, draw opacity={1.0}, line width={1}, solid}, colorbar={false}]
    \addplot[color={rgb,1:red,0.0;green,0.6056;blue,0.9787}, name path={356b3835-968e-4f22-9815-f7121540d081}, const plot, draw opacity={1.0}, line width={2}, dashed]
        table[row sep={\\}]
        {
            \\
            1.0  1.0  \\
            1.3333333333333317  0.1  \\
            2.3333333333333366  0.010000000000000002  \\
            2.3333333333333366  0.001  \\
            2.3333333333333366  0.0001  \\
            3.0  1.0e-5  \\
            3.333333333333332  1.0e-6  \\
            3.333333333333332  1.0e-7  \\
            3.333333333333332  1.0e-8  \\
            3.6666666666666696  1.0e-9  \\
            4.0  1.0e-10  \\
            4.0  1.0e-11  \\
            4.0  1.0e-12  \\
            4.0  1.0e-13  \\
            5.0  1.0e-14  \\
            5.0  1.0e-15  \\
            5.0  1.0e-16  \\
        }
        ;
    \addplot[color={rgb,1:red,0.0;green,0.0;blue,0.0}, name path={136abaa3-bcdd-48eb-9806-a47820dabb30}, draw opacity={1.0}, line width={1}, dotted]
        table[row sep={\\}]
        {
            \\
            6.333333333333333  1.0e-34  \\
            6.333333333333333  1.0e17  \\
        }
        ;
    \addplot[color={rgb,1:red,0.2422;green,0.6433;blue,0.3044}, name path={188e7053-7270-4edb-8c23-304848538f29}, const plot, draw opacity={1.0}, line width={2}, solid]
        table[row sep={\\}]
        {
            \\
            1.0  0.016554806653934297  \\
            1.3333333333333317  0.0007197002179759806  \\
            2.3333333333333366  5.3447982123582045e-8  \\
            2.3333333333333366  5.3447982123582045e-8  \\
            2.3333333333333366  5.3447982123582045e-8  \\
            3.0  1.5699008241418628e-9  \\
            3.333333333333332  9.724516856050285e-13  \\
            3.333333333333332  9.724516856050285e-13  \\
            3.333333333333332  9.724516856050285e-13  \\
            3.6666666666666696  1.5321400463081737e-14  \\
            4.0  4.967334462849976e-16  \\
            4.0  4.967334462849976e-16  \\
            4.0  4.967334462849976e-16  \\
            4.0  4.967334462849976e-16  \\
            5.0  4.0069831598718487e-16  \\
            5.0  4.0069831598718487e-16  \\
            5.0  6.303143724841514e-16  \\
        }
        ;
    \node[]  at (axis cs:6.27,1.0e-8) {$r_{9,9}$};
    \node[]  at (axis cs:1.0,0.08277403326967149) {$t_{2}$};
    \node[]  at (axis cs:1.3333333333333317,0.0035985010898799026) {$r_{2,1}$};
    \node[]  at (axis cs:2.3333333333333366,2.672399106179102e-7) {$r_{4,2}$};
    \node[]  at (axis cs:2.3333333333333366,2.672399106179102e-7) {$r_{4,2}$};
    \node[]  at (axis cs:2.3333333333333366,2.672399106179102e-7) {$r_{4,2}$};
    \node[]  at (axis cs:3.0,7.849504120709314e-9) {$t_{8}$};
    \node[]  at (axis cs:3.333333333333332,4.862258428025142e-12) {$r_{6,3}$};
    \node[]  at (axis cs:3.333333333333332,4.862258428025142e-12) {$r_{6,3}$};
    \node[]  at (axis cs:3.333333333333332,4.862258428025142e-12) {$r_{6,3}$};
    \node[]  at (axis cs:3.6666666666666696,7.660700231540869e-14) {$r_{6,4}$};
    \node[]  at (axis cs:4.0,2.4836672314249878e-15) {$t_{15}^{\left[16\right]}$};
    \node[]  at (axis cs:4.0,2.4836672314249878e-15) {$t_{15}^{\left[16\right]}$};
    \node[]  at (axis cs:4.0,2.4836672314249878e-15) {$t_{15}^{\left[16\right]}$};
    \node[]  at (axis cs:4.0,2.4836672314249878e-15) {$t_{15}^{\left[16\right]}$};
    \node[]  at (axis cs:5.0,2.0034915799359242e-15) {$t_{21}^{\left[24\right]}$};
    \node[]  at (axis cs:5.0,2.0034915799359242e-15) {$t_{21}^{\left[24\right]}$};
    \node[]  at (axis cs:5.0,3.1515718624207567e-14) {$t_{18}$};
\end{axis}
\end{tikzpicture}%
\end{subfigure}

\begin{subfigure}[b]{0.48\textwidth}
	\raggedleft
	\begin{tikzpicture}[/tikz/background rectangle/.style={fill={rgb,1:red,1.0;green,1.0;blue,1.0}, fill opacity={1.0}, draw opacity={1.0}}, show background rectangle]
\begin{axis}[point meta max={nan}, point meta min={nan}, legend cell align={left}, legend columns={1}, title={$\lVert hA \rVert_{1}=10.0$}, title style={at={{(0.5,1)}}, anchor={south}, color={rgb,1:red,0.0;green,0.0;blue,0.0}, draw opacity={1.0}, rotate={0.0}, align={center}}, legend style={color={rgb,1:red,0.0;green,0.0;blue,0.0}, draw opacity={1.0}, line width={1}, solid, fill={rgb,1:red,1.0;green,1.0;blue,1.0}, fill opacity={1.0}, text opacity={1.0}, text={rgb,1:red,0.0;green,0.0;blue,0.0}, cells={anchor={center}}, at={(1.02, 1)}, anchor={north west}}, axis background/.style={fill={rgb,1:red,1.0;green,1.0;blue,1.0}, opacity={1.0}}, anchor={north west}, scaled x ticks={false}, xlabel={cost,\ [$C$]}, x tick style={color={rgb,1:red,0.0;green,0.0;blue,0.0}, opacity={1.0}}, x tick label style={color={rgb,1:red,0.0;green,0.0;blue,0.0}, opacity={1.0}, rotate={0}}, xmajorgrids={true}, xmin={4.0}, xmax={8.453333333333333}, xticklabels={{$5$,$6$,$7$,$8$}}, xtick={{5.0,6.0,7.0,8.0}}, xtick align={inside}, xticklabel style={color={rgb,1:red,0.0;green,0.0;blue,0.0}, draw opacity={1.0}, rotate={0.0}}, x grid style={color={rgb,1:red,0.0;green,0.0;blue,0.0}, draw opacity={0.1}, line width={0.5}, solid}, axis x line*={left}, x axis line style={color={rgb,1:red,0.0;green,0.0;blue,0.0}, draw opacity={1.0}, line width={1}, solid}, scaled y ticks={false}, ylabel={$\lVert w_{\alpha}(h A)-e^{h A} \rVert_{1} / (\lVert hA \rVert_{1} \cdot \lVert e^{hA} \rVert_{1})$}, y tick style={color={rgb,1:red,0.0;green,0.0;blue,0.0}, opacity={1.0}}, y tick label style={color={rgb,1:red,0.0;green,0.0;blue,0.0}, opacity={1.0}, rotate={0}}, ymode={log}, log basis y={10}, ymajorgrids={true}, ymin={1.0e-17}, ymax={1.0}, yticklabels={{$10^{-16}$,$10^{-8}$,$10^{-12}$,$10^{-4}$}}, ytick={{1.0e-16,1.0e-8,1.0e-12,0.0001}}, ytick align={inside}, yticklabel style={color={rgb,1:red,0.0;green,0.0;blue,0.0}, draw opacity={1.0}, rotate={0.0}}, y grid style={color={rgb,1:red,0.0;green,0.0;blue,0.0}, draw opacity={0.1}, line width={0.5}, solid}, axis y line*={left}, y axis line style={color={rgb,1:red,0.0;green,0.0;blue,0.0}, draw opacity={1.0}, line width={1}, solid}, colorbar={false}]
    \addplot[color={rgb,1:red,0.0;green,0.6056;blue,0.9787}, name path={68a855a9-311e-431e-aff6-e9ecce01a406}, const plot, draw opacity={1.0}, line width={2}, dashed]
        table[row sep={\\}]
        {
            \\
            4.333333333333328  1.0  \\
            4.333333333333328  0.1  \\
            4.333333333333328  0.010000000000000002  \\
            5.333333333333327  0.001  \\
            5.333333333333327  0.0001  \\
            5.66666666666667  1.0e-5  \\
            6.0  1.0e-6  \\
            6.333333333333335  1.0e-7  \\
            6.666666666666664  1.0e-8  \\
            6.666666666666664  1.0e-9  \\
            6.666666666666664  1.0e-10  \\
            7.0  1.0e-11  \\
            7.66666666666666  1.0e-12  \\
            7.66666666666666  1.0e-13  \\
            7.66666666666666  1.0e-14  \\
            7.66666666666666  1.0e-15  \\
            7.66666666666666  1.0e-16  \\
        }
        ;
    \addplot[color={rgb,1:red,0.0;green,0.0;blue,0.0}, name path={eb0c91f7-000d-4bc5-b002-b156a375b7a1}, draw opacity={1.0}, line width={1}, dotted]
        table[row sep={\\}]
        {
            \\
            8.333333333333332  1.0e-34  \\
            8.333333333333332  1.0e17  \\
        }
        ;
    \addplot[color={rgb,1:red,0.2422;green,0.6433;blue,0.3044}, name path={6842c1b3-7bf7-47a5-8ee1-ec3c670e1934}, const plot, draw opacity={1.0}, line width={2}, solid]
        table[row sep={\\}]
        {
            \\
            4.333333333333328  6.84113224680917e-7  \\
            4.333333333333328  6.84113224680917e-7  \\
            4.333333333333328  6.125256218710389e-6  \\
            5.333333333333327  5.892944989287429e-10  \\
            5.333333333333327  1.6905770960941275e-9  \\
            5.66666666666667  7.256487865526352e-11  \\
            6.0  7.840411229215517e-15  \\
            6.333333333333335  1.8668084266689696e-13  \\
            6.666666666666664  1.2500164634245264e-13  \\
            6.666666666666664  1.2500164634245264e-13  \\
            6.666666666666664  1.2500164634245264e-13  \\
            7.0  1.448986007105137e-16  \\
            7.66666666666666  2.8939599911369614e-13  \\
            7.66666666666666  2.8939599911369614e-13  \\
            7.66666666666666  2.8939599911369614e-13  \\
            7.66666666666666  2.8939599911369614e-13  \\
            7.66666666666666  2.8939599911369614e-13  \\
        }
        ;
    \node[]  at (axis cs:8.249999999999998,1.0e-8) {$r_{13,13}$};
    \node[]  at (axis cs:4.333333333333328,3.4205661234045852e-6) {$r_{6,3}$};
    \node[]  at (axis cs:4.333333333333328,3.4205661234045852e-6) {$r_{6,3}$};
    \node[]  at (axis cs:4.333333333333328,3.0626281093551946e-5) {$r_{4,2}$};
    \node[]  at (axis cs:5.333333333333327,2.946472494643715e-9) {$r_{8,4}$};
    \node[]  at (axis cs:5.333333333333327,8.452885480470638e-9) {$r_{6,3}$};
    \node[]  at (axis cs:5.66666666666667,3.628243932763176e-10) {$r_{6,4}$};
    \node[]  at (axis cs:6.0,3.9202056146077585e-14) {$t_{15}^{\left[16\right]}$};
    \node[]  at (axis cs:6.333333333333335,9.334042133344847e-12) {$r_{8,4}$};
    \node[]  at (axis cs:6.666666666666664,6.250082317122632e-13) {$r_{12,8}$};
    \node[]  at (axis cs:6.666666666666664,6.250082317122632e-13) {$r_{12,8}$};
    \node[]  at (axis cs:6.666666666666664,6.250082317122632e-13) {$r_{12,8}$};
    \node[]  at (axis cs:7.0,7.244930035525685e-16) {$t_{21}^{\left[24\right]}$};
    \node[]  at (axis cs:7.66666666666666,1.4469799955684808e-12) {$r_{12,8}$};
    \node[]  at (axis cs:7.66666666666666,1.4469799955684808e-12) {$r_{12,8}$};
    \node[]  at (axis cs:7.66666666666666,1.4469799955684808e-12) {$r_{12,8}$};
    \node[]  at (axis cs:7.66666666666666,1.4469799955684808e-12) {$r_{12,8}$};
    \node[]  at (axis cs:7.66666666666666,1.4469799955684808e-12) {$r_{12,8}$};
\end{axis}
\end{tikzpicture}%
\end{subfigure}
~
\begin{subfigure}[b]{0.48\textwidth}
	\raggedright
	\begin{tikzpicture}[/tikz/background rectangle/.style={fill={rgb,1:red,1.0;green,1.0;blue,1.0}, fill opacity={1.0}, draw opacity={1.0}}, show background rectangle]
\begin{axis}[point meta max={nan}, point meta min={nan}, legend cell align={left}, legend columns={1}, title={$\lVert hA \rVert_{1}=100.0$}, title style={at={{(0.5,1)}}, anchor={south}, color={rgb,1:red,0.0;green,0.0;blue,0.0}, draw opacity={1.0}, rotate={0.0}, align={center}}, legend style={color={rgb,1:red,0.0;green,0.0;blue,0.0}, draw opacity={1.0}, line width={1}, solid, fill={rgb,1:red,1.0;green,1.0;blue,1.0}, fill opacity={1.0}, text opacity={1.0}, text={rgb,1:red,0.0;green,0.0;blue,0.0}, cells={anchor={center}}, at={(1.02, 1)}, anchor={north west}}, axis background/.style={fill={rgb,1:red,1.0;green,1.0;blue,1.0}, opacity={1.0}}, anchor={north west}, scaled x ticks={false}, xlabel={cost,\ [$C$]}, x tick style={color={rgb,1:red,0.0;green,0.0;blue,0.0}, opacity={1.0}}, x tick label style={color={rgb,1:red,0.0;green,0.0;blue,0.0}, opacity={1.0}, rotate={0}}, xmajorgrids={true}, xmin={7.0}, xmax={12.483333333333333}, xticklabels={{$8$,$9$,$10$,$11$,$12$}}, xtick={{8.0,9.0,10.0,11.0,12.0}}, xtick align={inside}, xticklabel style={color={rgb,1:red,0.0;green,0.0;blue,0.0}, draw opacity={1.0}, rotate={0.0}}, x grid style={color={rgb,1:red,0.0;green,0.0;blue,0.0}, draw opacity={0.1}, line width={0.5}, solid}, axis x line*={left}, x axis line style={color={rgb,1:red,0.0;green,0.0;blue,0.0}, draw opacity={1.0}, line width={1}, solid}, scaled y ticks={false}, ylabel={$\lVert w_{\alpha}(h A)-e^{h A} \rVert_{1} / (\lVert hA \rVert_{1} \cdot \lVert e^{hA} \rVert_{1})$}, y tick style={color={rgb,1:red,0.0;green,0.0;blue,0.0}, opacity={1.0}}, y tick label style={color={rgb,1:red,0.0;green,0.0;blue,0.0}, opacity={1.0}, rotate={0}}, ymode={log}, log basis y={10}, ymajorgrids={true}, ymin={1.0e-17}, ymax={1.0}, yticklabels={{$10^{-16}$,$10^{-8}$,$10^{-12}$,$10^{-4}$}}, ytick={{1.0e-16,1.0e-8,1.0e-12,0.0001}}, ytick align={inside}, yticklabel style={color={rgb,1:red,0.0;green,0.0;blue,0.0}, draw opacity={1.0}, rotate={0.0}}, y grid style={color={rgb,1:red,0.0;green,0.0;blue,0.0}, draw opacity={0.1}, line width={0.5}, solid}, axis y line*={left}, y axis line style={color={rgb,1:red,0.0;green,0.0;blue,0.0}, draw opacity={1.0}, line width={1}, solid}, colorbar={false}]
    \addplot[color={rgb,1:red,0.0;green,0.6056;blue,0.9787}, name path={7b45283e-d396-4d84-8948-2052fa88cb6a}, const plot, draw opacity={1.0}, line width={2}, dashed]
        table[row sep={\\}]
        {
            \\
            7.333333333333339  1.0  \\
            7.333333333333339  0.1  \\
            8.333333333333343  0.010000000000000002  \\
            8.333333333333343  0.001  \\
            8.666666666666655  0.0001  \\
            9.0  1.0e-5  \\
            9.333333333333346  1.0e-6  \\
            9.666666666666652  1.0e-7  \\
            9.666666666666652  1.0e-8  \\
            10.0  1.0e-9  \\
            10.0  1.0e-10  \\
            10.666666666666654  1.0e-11  \\
            10.666666666666654  1.0e-12  \\
            10.666666666666654  1.0e-13  \\
            10.666666666666654  1.0e-14  \\
            11.666666666666663  1.0e-15  \\
            11.666666666666663  1.0e-16  \\
        }
        ;
    \addplot[color={rgb,1:red,0.0;green,0.0;blue,0.0}, name path={16ce75c2-7a1b-4675-ac3e-bd396a8feced}, draw opacity={1.0}, line width={1}, dotted]
        table[row sep={\\}]
        {
            \\
            12.333333333333332  1.0e-34  \\
            12.333333333333332  1.0e17  \\
        }
        ;
    \addplot[color={rgb,1:red,0.2422;green,0.6433;blue,0.3044}, name path={a6b0221b-816a-4d8e-9037-d007611324aa}, const plot, draw opacity={1.0}, line width={2}, solid]
        table[row sep={\\}]
        {
            \\
            7.333333333333339  1.922920159002151e-5  \\
            7.333333333333339  1.922920159002151e-5  \\
            8.333333333333343  1.0166835785144016e-8  \\
            8.333333333333343  1.0166835785144016e-8  \\
            8.666666666666655  5.620786367074017e-10  \\
            9.0  1.5172483222416145e-13  \\
            9.333333333333346  2.1870926187856914e-12  \\
            9.666666666666652  6.452789572185826e-14  \\
            9.666666666666652  6.452789572185826e-14  \\
            10.0  7.53383397801844e-17  \\
            10.0  1.2596676701773128e-16  \\
            10.666666666666654  1.835142446537144e-13  \\
            10.666666666666654  1.835142446537144e-13  \\
            10.666666666666654  1.835142446537144e-13  \\
            10.666666666666654  1.835142446537144e-13  \\
            11.666666666666663  3.258408514302809e-13  \\
            11.666666666666663  3.258408514302809e-13  \\
        }
        ;
    \node[]  at (axis cs:12.209999999999999,1.0e-8) {$r_{13,13}$};
    \node[]  at (axis cs:7.333333333333339,9.614600795010754e-5) {$r_{4,2}$};
    \node[]  at (axis cs:7.333333333333339,9.614600795010754e-5) {$r_{4,2}$};
    \node[]  at (axis cs:8.333333333333343,5.0834178925720086e-8) {$r_{6,3}$};
    \node[]  at (axis cs:8.333333333333343,5.0834178925720086e-8) {$r_{6,3}$};
    \node[]  at (axis cs:8.666666666666655,2.8103931835370085e-9) {$r_{6,4}$};
    \node[]  at (axis cs:9.0,7.586241611208072e-13) {$t_{15}^{\left[16\right]}$};
    \node[]  at (axis cs:9.333333333333346,1.0935463093928457e-11) {$r_{8,4}$};
    \node[]  at (axis cs:9.666666666666652,3.226394786092913e-13) {$r_{12,8}$};
    \node[]  at (axis cs:9.666666666666652,3.226394786092913e-13) {$r_{12,8}$};
    \node[]  at (axis cs:10.0,0.76691698900922e-16) {$t_{21}^{\left[24\right]}$};
    \node[]  at (axis cs:10.0,6.298338350886564e-15) {$t_{15}^{\left[16\right]}$};
    \node[]  at (axis cs:10.666666666666654,9.17571223268572e-13) {$r_{12,8}$};
    \node[]  at (axis cs:10.666666666666654,9.17571223268572e-13) {$r_{12,8}$};
    \node[]  at (axis cs:10.666666666666654,9.17571223268572e-13) {$r_{12,8}$};
    \node[]  at (axis cs:10.666666666666654,9.17571223268572e-13) {$r_{12,8}$};
    \node[]  at (axis cs:11.666666666666663,1.6292042571514045e-12) {$r_{12,8}$};
    \node[]  at (axis cs:11.666666666666663,1.6292042571514045e-12) {$r_{12,8}$};
\end{axis}
\end{tikzpicture}%
\end{subfigure}
\caption{Comparison of tolerance (dashed) and practical average cost plotted vs. total computational cost of superdiagonal methods.
	The thin vertical dotted line represents the cost in the current software that implements \cite{higham09tsa}.
	Note that the distance between the dotted vertical line and the solid one represents cost savings.}
\label{fig:error_vs_cost_superdiag}
\end{figure}

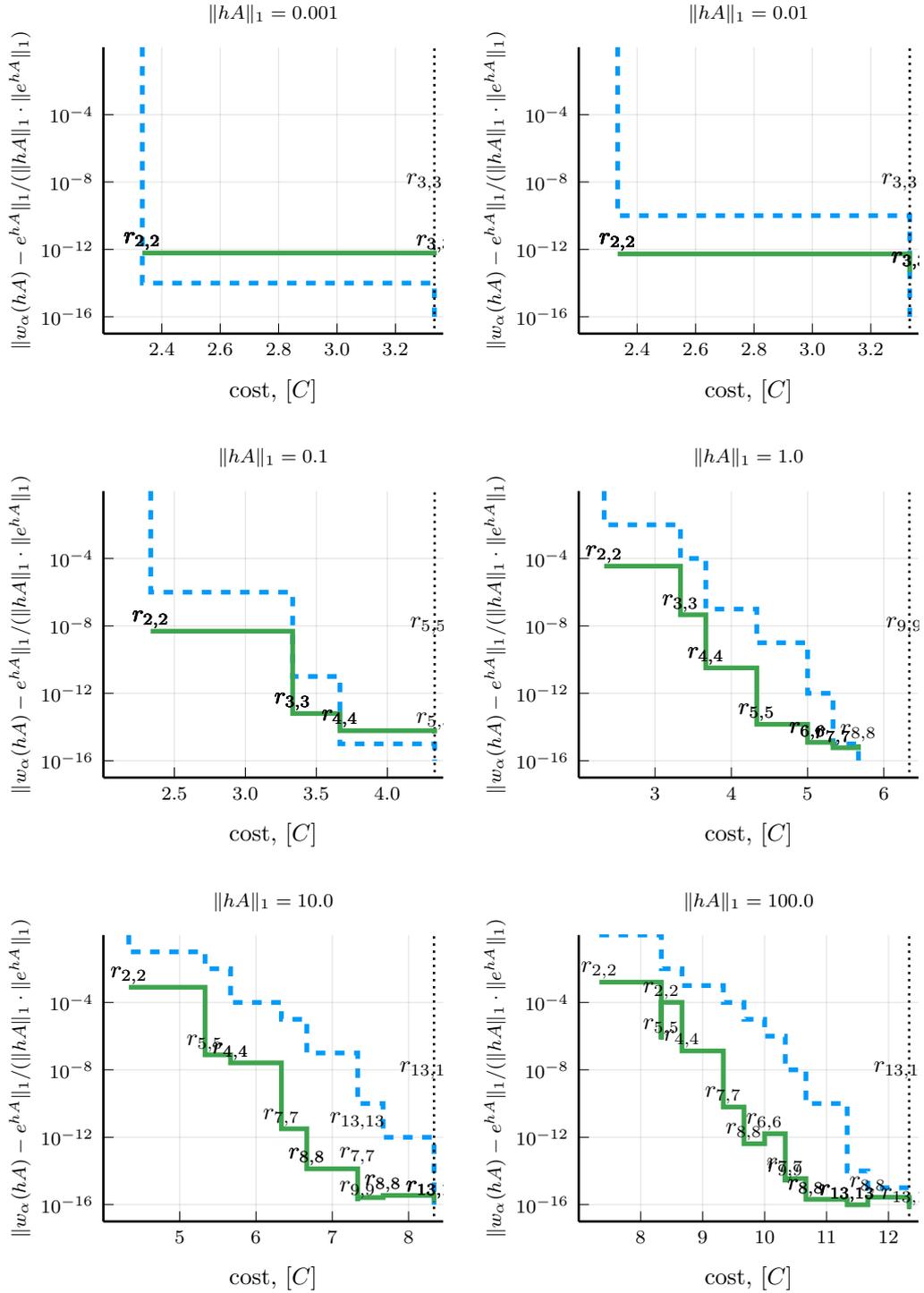
\begin{figure}[p]
	\begin{subfigure}[b]{0.48\textwidth}
		\raggedleft
		\begin{tikzpicture}[/tikz/background rectangle/.style={fill={rgb,1:red,1.0;green,1.0;blue,1.0}, fill opacity={1.0}, draw opacity={1.0}}, show background rectangle]
\begin{axis}[point meta max={nan}, point meta min={nan}, legend cell align={left}, legend columns={1}, title={$\lVert hA \rVert_{1}=0.001$}, title style={at={{(0.5,1)}}, anchor={south}, color={rgb,1:red,0.0;green,0.0;blue,0.0}, draw opacity={1.0}, rotate={0.0}, align={center}}, legend style={color={rgb,1:red,0.0;green,0.0;blue,0.0}, draw opacity={1.0}, line width={1}, solid, fill={rgb,1:red,1.0;green,1.0;blue,1.0}, fill opacity={1.0}, text opacity={1.0}, text={rgb,1:red,0.0;green,0.0;blue,0.0}, cells={anchor={center}}, at={(1.02, 1)}, anchor={north west}}, axis background/.style={fill={rgb,1:red,1.0;green,1.0;blue,1.0}, opacity={1.0}}, anchor={north west}, scaled x ticks={false}, xlabel={cost,\ [$C$]}, x tick style={color={rgb,1:red,0.0;green,0.0;blue,0.0}, opacity={1.0}}, x tick label style={color={rgb,1:red,0.0;green,0.0;blue,0.0}, opacity={1.0}, rotate={0}}, xmajorgrids={true}, xmin={2.2}, xmax={3.3633333333333333}, xticklabels={{$2.4$,$2.6$,$2.8$,$3.0$,$3.2$}}, xtick={{2.4000000000000004,2.6,2.8000000000000003,3.0,3.2}}, xtick align={inside}, xticklabel style={color={rgb,1:red,0.0;green,0.0;blue,0.0}, draw opacity={1.0}, rotate={0.0}}, x grid style={color={rgb,1:red,0.0;green,0.0;blue,0.0}, draw opacity={0.1}, line width={0.5}, solid}, axis x line*={left}, x axis line style={color={rgb,1:red,0.0;green,0.0;blue,0.0}, draw opacity={1.0}, line width={1}, solid}, scaled y ticks={false}, ylabel={$\lVert w_{\alpha}(h A)-e^{h A} \rVert_{1} / (\lVert hA \rVert_{1} \cdot \lVert e^{hA} \rVert_{1})$}, y tick style={color={rgb,1:red,0.0;green,0.0;blue,0.0}, opacity={1.0}}, y tick label style={color={rgb,1:red,0.0;green,0.0;blue,0.0}, opacity={1.0}, rotate={0}}, ymode={log}, log basis y={10}, ymajorgrids={true}, ymin={1.0e-17}, ymax={1.0}, yticklabels={{$10^{-16}$,$10^{-8}$,$10^{-12}$,$10^{-4}$}}, ytick={{1.0e-16,1.0e-8,1.0e-12,0.0001}}, ytick align={inside}, yticklabel style={color={rgb,1:red,0.0;green,0.0;blue,0.0}, draw opacity={1.0}, rotate={0.0}}, y grid style={color={rgb,1:red,0.0;green,0.0;blue,0.0}, draw opacity={0.1}, line width={0.5}, solid}, axis y line*={left}, y axis line style={color={rgb,1:red,0.0;green,0.0;blue,0.0}, draw opacity={1.0}, line width={1}, solid}, colorbar={false}]
    \addplot[color={rgb,1:red,0.0;green,0.6056;blue,0.9787}, name path={dcf1c608-e192-48b7-8730-b9dae05bf7a0}, const plot, draw opacity={1.0}, line width={2}, dashed]
        table[row sep={\\}]
        {
            \\
            2.3333333333333366  1.0  \\
            2.3333333333333366  0.1  \\
            2.3333333333333366  0.010000000000000002  \\
            2.3333333333333366  0.001  \\
            2.3333333333333366  0.0001  \\
            2.3333333333333366  1.0e-5  \\
            2.3333333333333366  1.0e-6  \\
            2.3333333333333366  1.0e-7  \\
            2.3333333333333366  1.0e-8  \\
            2.3333333333333366  1.0e-9  \\
            2.3333333333333366  1.0e-10  \\
            2.3333333333333366  1.0e-11  \\
            2.3333333333333366  1.0e-12  \\
            2.3333333333333366  1.0e-13  \\
            2.3333333333333366  1.0e-14  \\
            3.333333333333332  1.0e-15  \\
            3.333333333333332  1.0e-16  \\
        }
        ;
    \addplot[color={rgb,1:red,0.0;green,0.0;blue,0.0}, name path={6471fc45-29a5-4d74-9450-f7d64f039104}, draw opacity={1.0}, line width={1}, dotted]
        table[row sep={\\}]
        {
            \\
            3.3333333333333335  1.0e-34  \\
            3.3333333333333335  1.0e17  \\
        }
        ;
    \addplot[color={rgb,1:red,0.2422;green,0.6433;blue,0.3044}, name path={12fae452-dc64-4f43-92a4-ffed89e980db}, const plot, draw opacity={1.0}, line width={2}, solid]
        table[row sep={\\}]
        {
            \\
            2.3333333333333366  6.102437563442668e-13  \\
            2.3333333333333366  6.102437563442668e-13  \\
            2.3333333333333366  6.102437563442668e-13  \\
            2.3333333333333366  6.102437563442668e-13  \\
            2.3333333333333366  6.102437563442668e-13  \\
            2.3333333333333366  6.102437563442668e-13  \\
            2.3333333333333366  6.102437563442668e-13  \\
            2.3333333333333366  6.102437563442668e-13  \\
            2.3333333333333366  6.102437563442668e-13  \\
            2.3333333333333366  6.102437563442668e-13  \\
            2.3333333333333366  6.102437563442668e-13  \\
            2.3333333333333366  6.102437563442668e-13  \\
            2.3333333333333366  6.102437563442668e-13  \\
            2.3333333333333366  6.102437563442668e-13  \\
            2.3333333333333366  6.102437563442668e-13  \\
            3.333333333333332  4.4159194241260415e-13  \\
            3.333333333333332  4.4159194241260415e-13  \\
        }
        ;
    \node[]  at (axis cs:3.3000000000000003,1.0e-8) {$r_{3,3}$};
    \node[]  at (axis cs:2.3333333333333366,3.0512187817213337e-12) {$r_{2,2}$};
    \node[]  at (axis cs:2.3333333333333366,3.0512187817213337e-12) {$r_{2,2}$};
    \node[]  at (axis cs:2.3333333333333366,3.0512187817213337e-12) {$r_{2,2}$};
    \node[]  at (axis cs:2.3333333333333366,3.0512187817213337e-12) {$r_{2,2}$};
    \node[]  at (axis cs:2.3333333333333366,3.0512187817213337e-12) {$r_{2,2}$};
    \node[]  at (axis cs:2.3333333333333366,3.0512187817213337e-12) {$r_{2,2}$};
    \node[]  at (axis cs:2.3333333333333366,3.0512187817213337e-12) {$r_{2,2}$};
    \node[]  at (axis cs:2.3333333333333366,3.0512187817213337e-12) {$r_{2,2}$};
    \node[]  at (axis cs:2.3333333333333366,3.0512187817213337e-12) {$r_{2,2}$};
    \node[]  at (axis cs:2.3333333333333366,3.0512187817213337e-12) {$r_{2,2}$};
    \node[]  at (axis cs:2.3333333333333366,3.0512187817213337e-12) {$r_{2,2}$};
    \node[]  at (axis cs:2.3333333333333366,3.0512187817213337e-12) {$r_{2,2}$};
    \node[]  at (axis cs:2.3333333333333366,3.0512187817213337e-12) {$r_{2,2}$};
    \node[]  at (axis cs:2.3333333333333366,3.0512187817213337e-12) {$r_{2,2}$};
    \node[]  at (axis cs:2.3333333333333366,3.0512187817213337e-12) {$r_{2,2}$};
    \node[]  at (axis cs:3.333333333333332,2.207959712063021e-12) {$r_{3,3}$};
    \node[]  at (axis cs:3.333333333333332,2.207959712063021e-12) {$r_{3,3}$};
\end{axis}
\end{tikzpicture}%
	\end{subfigure}
	~
	\begin{subfigure}[b]{0.48\textwidth}
		\raggedright
		\begin{tikzpicture}[/tikz/background rectangle/.style={fill={rgb,1:red,1.0;green,1.0;blue,1.0}, fill opacity={1.0}, draw opacity={1.0}}, show background rectangle]
\begin{axis}[point meta max={nan}, point meta min={nan}, legend cell align={left}, legend columns={1}, title={$\lVert hA \rVert_{1}=0.01$}, title style={at={{(0.5,1)}}, anchor={south}, color={rgb,1:red,0.0;green,0.0;blue,0.0}, draw opacity={1.0}, rotate={0.0}, align={center}}, legend style={color={rgb,1:red,0.0;green,0.0;blue,0.0}, draw opacity={1.0}, line width={1}, solid, fill={rgb,1:red,1.0;green,1.0;blue,1.0}, fill opacity={1.0}, text opacity={1.0}, text={rgb,1:red,0.0;green,0.0;blue,0.0}, cells={anchor={center}}, at={(1.02, 1)}, anchor={north west}}, axis background/.style={fill={rgb,1:red,1.0;green,1.0;blue,1.0}, opacity={1.0}}, anchor={north west}, scaled x ticks={false}, xlabel={cost,\ [$C$]}, x tick style={color={rgb,1:red,0.0;green,0.0;blue,0.0}, opacity={1.0}}, x tick label style={color={rgb,1:red,0.0;green,0.0;blue,0.0}, opacity={1.0}, rotate={0}}, xmajorgrids={true}, xmin={2.2}, xmax={3.3633333333333333}, xticklabels={{$2.4$,$2.6$,$2.8$,$3.0$,$3.2$}}, xtick={{2.4000000000000004,2.6,2.8000000000000003,3.0,3.2}}, xtick align={inside}, xticklabel style={color={rgb,1:red,0.0;green,0.0;blue,0.0}, draw opacity={1.0}, rotate={0.0}}, x grid style={color={rgb,1:red,0.0;green,0.0;blue,0.0}, draw opacity={0.1}, line width={0.5}, solid}, axis x line*={left}, x axis line style={color={rgb,1:red,0.0;green,0.0;blue,0.0}, draw opacity={1.0}, line width={1}, solid}, scaled y ticks={false}, ylabel={$\lVert w_{\alpha}(h A)-e^{h A} \rVert_{1} / (\lVert hA \rVert_{1} \cdot \lVert e^{hA} \rVert_{1})$}, y tick style={color={rgb,1:red,0.0;green,0.0;blue,0.0}, opacity={1.0}}, y tick label style={color={rgb,1:red,0.0;green,0.0;blue,0.0}, opacity={1.0}, rotate={0}}, ymode={log}, log basis y={10}, ymajorgrids={true}, ymin={1.0e-17}, ymax={1.0}, yticklabels={{$10^{-16}$,$10^{-8}$,$10^{-12}$,$10^{-4}$}}, ytick={{1.0e-16,1.0e-8,1.0e-12,0.0001}}, ytick align={inside}, yticklabel style={color={rgb,1:red,0.0;green,0.0;blue,0.0}, draw opacity={1.0}, rotate={0.0}}, y grid style={color={rgb,1:red,0.0;green,0.0;blue,0.0}, draw opacity={0.1}, line width={0.5}, solid}, axis y line*={left}, y axis line style={color={rgb,1:red,0.0;green,0.0;blue,0.0}, draw opacity={1.0}, line width={1}, solid}, colorbar={false}]
    \addplot[color={rgb,1:red,0.0;green,0.6056;blue,0.9787}, name path={83617efa-62f8-4342-8c63-94e26e8ccb06}, const plot, draw opacity={1.0}, line width={2}, dashed]
        table[row sep={\\}]
        {
            \\
            2.3333333333333366  1.0  \\
            2.3333333333333366  0.1  \\
            2.3333333333333366  0.010000000000000002  \\
            2.3333333333333366  0.001  \\
            2.3333333333333366  0.0001  \\
            2.3333333333333366  1.0e-5  \\
            2.3333333333333366  1.0e-6  \\
            2.3333333333333366  1.0e-7  \\
            2.3333333333333366  1.0e-8  \\
            2.3333333333333366  1.0e-9  \\
            2.3333333333333366  1.0e-10  \\
            3.333333333333332  1.0e-11  \\
            3.333333333333332  1.0e-12  \\
            3.333333333333332  1.0e-13  \\
            3.333333333333332  1.0e-14  \\
            3.333333333333332  1.0e-15  \\
            3.333333333333332  1.0e-16  \\
        }
        ;
    \addplot[color={rgb,1:red,0.0;green,0.0;blue,0.0}, name path={4fce6010-3cdd-4131-8dd9-d8158ec03bd5}, draw opacity={1.0}, line width={1}, dotted]
        table[row sep={\\}]
        {
            \\
            3.3333333333333335  1.0e-34  \\
            3.3333333333333335  1.0e17  \\
        }
        ;
    \addplot[color={rgb,1:red,0.2422;green,0.6433;blue,0.3044}, name path={ebd15634-26d7-45de-bacb-dad9b01ccfcc}, const plot, draw opacity={1.0}, line width={2}, solid]
        table[row sep={\\}]
        {
            \\
            2.3333333333333366  5.388322267276272e-13  \\
            2.3333333333333366  5.388322267276272e-13  \\
            2.3333333333333366  5.388322267276272e-13  \\
            2.3333333333333366  5.388322267276272e-13  \\
            2.3333333333333366  5.388322267276272e-13  \\
            2.3333333333333366  5.388322267276272e-13  \\
            2.3333333333333366  5.388322267276272e-13  \\
            2.3333333333333366  5.388322267276272e-13  \\
            2.3333333333333366  5.388322267276272e-13  \\
            2.3333333333333366  5.388322267276272e-13  \\
            2.3333333333333366  5.388322267276272e-13  \\
            3.333333333333332  4.380562585859002e-14  \\
            3.333333333333332  4.380562585859002e-14  \\
            3.333333333333332  4.380562585859002e-14  \\
            3.333333333333332  4.380562585859002e-14  \\
            3.333333333333332  4.380562585859002e-14  \\
            3.333333333333332  4.380562585859002e-14  \\
        }
        ;
    \node[]  at (axis cs:3.3000000000000003,1.0e-8) {$r_{3,3}$};
    \node[]  at (axis cs:2.3333333333333366,2.6941611336381362e-12) {$r_{2,2}$};
    \node[]  at (axis cs:2.3333333333333366,2.6941611336381362e-12) {$r_{2,2}$};
    \node[]  at (axis cs:2.3333333333333366,2.6941611336381362e-12) {$r_{2,2}$};
    \node[]  at (axis cs:2.3333333333333366,2.6941611336381362e-12) {$r_{2,2}$};
    \node[]  at (axis cs:2.3333333333333366,2.6941611336381362e-12) {$r_{2,2}$};
    \node[]  at (axis cs:2.3333333333333366,2.6941611336381362e-12) {$r_{2,2}$};
    \node[]  at (axis cs:2.3333333333333366,2.6941611336381362e-12) {$r_{2,2}$};
    \node[]  at (axis cs:2.3333333333333366,2.6941611336381362e-12) {$r_{2,2}$};
    \node[]  at (axis cs:2.3333333333333366,2.6941611336381362e-12) {$r_{2,2}$};
    \node[]  at (axis cs:2.3333333333333366,2.6941611336381362e-12) {$r_{2,2}$};
    \node[]  at (axis cs:2.3333333333333366,2.6941611336381362e-12) {$r_{2,2}$};
    \node[]  at (axis cs:3.333333333333332,2.190281292929501e-13) {$r_{3,3}$};
    \node[]  at (axis cs:3.333333333333332,2.190281292929501e-13) {$r_{3,3}$};
    \node[]  at (axis cs:3.333333333333332,2.190281292929501e-13) {$r_{3,3}$};
    \node[]  at (axis cs:3.333333333333332,2.190281292929501e-13) {$r_{3,3}$};
    \node[]  at (axis cs:3.333333333333332,2.190281292929501e-13) {$r_{3,3}$};
    \node[]  at (axis cs:3.333333333333332,2.190281292929501e-13) {$r_{3,3}$};
\end{axis}
\end{tikzpicture}%
	\end{subfigure}

	\begin{subfigure}[b]{0.48\textwidth}
		\raggedleft
		\begin{tikzpicture}[/tikz/background rectangle/.style={fill={rgb,1:red,1.0;green,1.0;blue,1.0}, fill opacity={1.0}, draw opacity={1.0}}, show background rectangle]
\begin{axis}[point meta max={nan}, point meta min={nan}, legend cell align={left}, legend columns={1}, title={$\lVert hA \rVert_{1}=0.1$}, title style={at={{(0.5,1)}}, anchor={south}, color={rgb,1:red,0.0;green,0.0;blue,0.0}, draw opacity={1.0}, rotate={0.0}, align={center}}, legend style={color={rgb,1:red,0.0;green,0.0;blue,0.0}, draw opacity={1.0}, line width={1}, solid, fill={rgb,1:red,1.0;green,1.0;blue,1.0}, fill opacity={1.0}, text opacity={1.0}, text={rgb,1:red,0.0;green,0.0;blue,0.0}, cells={anchor={center}}, at={(1.02, 1)}, anchor={north west}}, axis background/.style={fill={rgb,1:red,1.0;green,1.0;blue,1.0}, opacity={1.0}}, anchor={north west}, scaled x ticks={false}, xlabel={cost,\ [$C$]}, x tick style={color={rgb,1:red,0.0;green,0.0;blue,0.0}, opacity={1.0}}, x tick label style={color={rgb,1:red,0.0;green,0.0;blue,0.0}, opacity={1.0}, rotate={0}}, xmajorgrids={true}, xmin={2.0}, xmax={4.393333333333333}, xticklabels={{$2.5$,$3.0$,$3.5$,$4.0$}}, xtick={{2.5,3.0,3.5,4.0}}, xtick align={inside}, xticklabel style={color={rgb,1:red,0.0;green,0.0;blue,0.0}, draw opacity={1.0}, rotate={0.0}}, x grid style={color={rgb,1:red,0.0;green,0.0;blue,0.0}, draw opacity={0.1}, line width={0.5}, solid}, axis x line*={left}, x axis line style={color={rgb,1:red,0.0;green,0.0;blue,0.0}, draw opacity={1.0}, line width={1}, solid}, scaled y ticks={false}, ylabel={$\lVert w_{\alpha}(h A)-e^{h A} \rVert_{1} / (\lVert hA \rVert_{1} \cdot \lVert e^{hA} \rVert_{1})$}, y tick style={color={rgb,1:red,0.0;green,0.0;blue,0.0}, opacity={1.0}}, y tick label style={color={rgb,1:red,0.0;green,0.0;blue,0.0}, opacity={1.0}, rotate={0}}, ymode={log}, log basis y={10}, ymajorgrids={true}, ymin={1.0e-17}, ymax={1.0}, yticklabels={{$10^{-16}$,$10^{-8}$,$10^{-12}$,$10^{-4}$}}, ytick={{1.0e-16,1.0e-8,1.0e-12,0.0001}}, ytick align={inside}, yticklabel style={color={rgb,1:red,0.0;green,0.0;blue,0.0}, draw opacity={1.0}, rotate={0.0}}, y grid style={color={rgb,1:red,0.0;green,0.0;blue,0.0}, draw opacity={0.1}, line width={0.5}, solid}, axis y line*={left}, y axis line style={color={rgb,1:red,0.0;green,0.0;blue,0.0}, draw opacity={1.0}, line width={1}, solid}, colorbar={false}]
    \addplot[color={rgb,1:red,0.0;green,0.6056;blue,0.9787}, name path={e71a7728-3dd5-423c-aee9-01328421db3f}, const plot, draw opacity={1.0}, line width={2}, dashed]
        table[row sep={\\}]
        {
            \\
            2.3333333333333366  1.0  \\
            2.3333333333333366  0.1  \\
            2.3333333333333366  0.010000000000000002  \\
            2.3333333333333366  0.001  \\
            2.3333333333333366  0.0001  \\
            2.3333333333333366  1.0e-5  \\
            2.3333333333333366  1.0e-6  \\
            3.333333333333332  1.0e-7  \\
            3.333333333333332  1.0e-8  \\
            3.333333333333332  1.0e-9  \\
            3.333333333333332  1.0e-10  \\
            3.333333333333332  1.0e-11  \\
            3.6666666666666696  1.0e-12  \\
            3.6666666666666696  1.0e-13  \\
            3.6666666666666696  1.0e-14  \\
            3.6666666666666696  1.0e-15  \\
            4.333333333333328  1.0e-16  \\
        }
        ;
    \addplot[color={rgb,1:red,0.0;green,0.0;blue,0.0}, name path={da104181-9423-4d07-8659-078eedb0489c}, draw opacity={1.0}, line width={1}, dotted]
        table[row sep={\\}]
        {
            \\
            4.333333333333333  1.0e-34  \\
            4.333333333333333  1.0e17  \\
        }
        ;
    \addplot[color={rgb,1:red,0.2422;green,0.6433;blue,0.3044}, name path={708237ac-0061-4c4e-8492-f4ee4283865c}, const plot, draw opacity={1.0}, line width={2}, solid]
        table[row sep={\\}]
        {
            \\
            2.3333333333333366  4.856919845872918e-9  \\
            2.3333333333333366  4.856919845872918e-9  \\
            2.3333333333333366  4.856919845872918e-9  \\
            2.3333333333333366  4.856919845872918e-9  \\
            2.3333333333333366  4.856919845872918e-9  \\
            2.3333333333333366  4.856919845872918e-9  \\
            2.3333333333333366  4.856919845872918e-9  \\
            3.333333333333332  6.31775897314289e-14  \\
            3.333333333333332  6.31775897314289e-14  \\
            3.333333333333332  6.31775897314289e-14  \\
            3.333333333333332  6.31775897314289e-14  \\
            3.333333333333332  6.31775897314289e-14  \\
            3.6666666666666696  6.0398069061415034e-15  \\
            3.6666666666666696  6.0398069061415034e-15  \\
            3.6666666666666696  6.0398069061415034e-15  \\
            3.6666666666666696  6.0398069061415034e-15  \\
            4.333333333333328  4.189708149646262e-15  \\
        }
        ;
    \node[]  at (axis cs:4.29,1.0e-8) {$r_{5,5}$};
    \node[]  at (axis cs:2.3333333333333366,2.428459922936459e-8) {$r_{2,2}$};
    \node[]  at (axis cs:2.3333333333333366,2.428459922936459e-8) {$r_{2,2}$};
    \node[]  at (axis cs:2.3333333333333366,2.428459922936459e-8) {$r_{2,2}$};
    \node[]  at (axis cs:2.3333333333333366,2.428459922936459e-8) {$r_{2,2}$};
    \node[]  at (axis cs:2.3333333333333366,2.428459922936459e-8) {$r_{2,2}$};
    \node[]  at (axis cs:2.3333333333333366,2.428459922936459e-8) {$r_{2,2}$};
    \node[]  at (axis cs:2.3333333333333366,2.428459922936459e-8) {$r_{2,2}$};
    \node[]  at (axis cs:3.333333333333332,3.1588794865714453e-13) {$r_{3,3}$};
    \node[]  at (axis cs:3.333333333333332,3.1588794865714453e-13) {$r_{3,3}$};
    \node[]  at (axis cs:3.333333333333332,3.1588794865714453e-13) {$r_{3,3}$};
    \node[]  at (axis cs:3.333333333333332,3.1588794865714453e-13) {$r_{3,3}$};
    \node[]  at (axis cs:3.333333333333332,3.1588794865714453e-13) {$r_{3,3}$};
    \node[]  at (axis cs:3.6666666666666696,3.0199034530707515e-14) {$r_{4,4}$};
    \node[]  at (axis cs:3.6666666666666696,3.0199034530707515e-14) {$r_{4,4}$};
    \node[]  at (axis cs:3.6666666666666696,3.0199034530707515e-14) {$r_{4,4}$};
    \node[]  at (axis cs:3.6666666666666696,3.0199034530707515e-14) {$r_{4,4}$};
    \node[]  at (axis cs:4.333333333333328,2.094854074823131e-14) {$r_{5,5}$};
\end{axis}
\end{tikzpicture}%
	\end{subfigure}
	~
	\begin{subfigure}[b]{0.48\textwidth}
		\raggedright
		\begin{tikzpicture}[/tikz/background rectangle/.style={fill={rgb,1:red,1.0;green,1.0;blue,1.0}, fill opacity={1.0}, draw opacity={1.0}}, show background rectangle]
\begin{axis}[point meta max={nan}, point meta min={nan}, legend cell align={left}, legend columns={1}, title={$\lVert hA \rVert_{1}=1.0$}, title style={at={{(0.5,1)}}, anchor={south}, color={rgb,1:red,0.0;green,0.0;blue,0.0}, draw opacity={1.0}, rotate={0.0}, align={center}}, legend style={color={rgb,1:red,0.0;green,0.0;blue,0.0}, draw opacity={1.0}, line width={1}, solid, fill={rgb,1:red,1.0;green,1.0;blue,1.0}, fill opacity={1.0}, text opacity={1.0}, text={rgb,1:red,0.0;green,0.0;blue,0.0}, cells={anchor={center}}, at={(1.02, 1)}, anchor={north west}}, axis background/.style={fill={rgb,1:red,1.0;green,1.0;blue,1.0}, opacity={1.0}}, anchor={north west}, scaled x ticks={false}, xlabel={cost,\ [$C$]}, x tick style={color={rgb,1:red,0.0;green,0.0;blue,0.0}, opacity={1.0}}, x tick label style={color={rgb,1:red,0.0;green,0.0;blue,0.0}, opacity={1.0}, rotate={0}}, xmajorgrids={true}, xmin={2.0}, xmax={6.453333333333333}, xticklabels={{$3$,$4$,$5$,$6$}}, xtick={{3.0,4.0,5.0,6.0}}, xtick align={inside}, xticklabel style={color={rgb,1:red,0.0;green,0.0;blue,0.0}, draw opacity={1.0}, rotate={0.0}}, x grid style={color={rgb,1:red,0.0;green,0.0;blue,0.0}, draw opacity={0.1}, line width={0.5}, solid}, axis x line*={left}, x axis line style={color={rgb,1:red,0.0;green,0.0;blue,0.0}, draw opacity={1.0}, line width={1}, solid}, scaled y ticks={false}, ylabel={$\lVert w_{\alpha}(h A)-e^{h A} \rVert_{1} / (\lVert hA \rVert_{1} \cdot \lVert e^{hA} \rVert_{1})$}, y tick style={color={rgb,1:red,0.0;green,0.0;blue,0.0}, opacity={1.0}}, y tick label style={color={rgb,1:red,0.0;green,0.0;blue,0.0}, opacity={1.0}, rotate={0}}, ymode={log}, log basis y={10}, ymajorgrids={true}, ymin={1.0e-17}, ymax={1.0}, yticklabels={{$10^{-16}$,$10^{-8}$,$10^{-12}$,$10^{-4}$}}, ytick={{1.0e-16,1.0e-8,1.0e-12,0.0001}}, ytick align={inside}, yticklabel style={color={rgb,1:red,0.0;green,0.0;blue,0.0}, draw opacity={1.0}, rotate={0.0}}, y grid style={color={rgb,1:red,0.0;green,0.0;blue,0.0}, draw opacity={0.1}, line width={0.5}, solid}, axis y line*={left}, y axis line style={color={rgb,1:red,0.0;green,0.0;blue,0.0}, draw opacity={1.0}, line width={1}, solid}, colorbar={false}]
    \addplot[color={rgb,1:red,0.0;green,0.6056;blue,0.9787}, name path={edc7250c-9f42-45e9-aed3-fc117f948b8f}, const plot, draw opacity={1.0}, line width={2}, dashed]
        table[row sep={\\}]
        {
            \\
            2.3333333333333366  1.0  \\
            2.3333333333333366  0.1  \\
            2.3333333333333366  0.010000000000000002  \\
            3.333333333333332  0.001  \\
            3.333333333333332  0.0001  \\
            3.6666666666666696  1.0e-5  \\
            3.6666666666666696  1.0e-6  \\
            3.6666666666666696  1.0e-7  \\
            4.333333333333328  1.0e-8  \\
            4.333333333333328  1.0e-9  \\
            5.0  1.0e-10  \\
            5.0  1.0e-11  \\
            5.0  1.0e-12  \\
            5.333333333333327  1.0e-13  \\
            5.333333333333327  1.0e-14  \\
            5.333333333333327  1.0e-15  \\
            5.66666666666667  1.0e-16  \\
        }
        ;
    \addplot[color={rgb,1:red,0.0;green,0.0;blue,0.0}, name path={ed1db9c2-fa7f-4e36-b599-31f52e4bdc36}, draw opacity={1.0}, line width={1}, dotted]
        table[row sep={\\}]
        {
            \\
            6.333333333333333  1.0e-34  \\
            6.333333333333333  1.0e17  \\
        }
        ;
    \addplot[color={rgb,1:red,0.2422;green,0.6433;blue,0.3044}, name path={e2224fe9-fd21-441c-991f-d12c552d324e}, const plot, draw opacity={1.0}, line width={2}, solid]
        table[row sep={\\}]
        {
            \\
            2.3333333333333366  3.548837479725014e-5  \\
            2.3333333333333366  3.548837479725014e-5  \\
            2.3333333333333366  3.548837479725014e-5  \\
            3.333333333333332  4.5669926928095953e-8  \\
            3.333333333333332  4.5669926928095953e-8  \\
            3.6666666666666696  3.228885776428477e-11  \\
            3.6666666666666696  3.228885776428477e-11  \\
            3.6666666666666696  3.228885776428477e-11  \\
            4.333333333333328  1.4458920107201225e-14  \\
            4.333333333333328  1.4458920107201225e-14  \\
            5.0  1.229701245866498e-15  \\
            5.0  1.229701245866498e-15  \\
            5.0  1.229701245866498e-15  \\
            5.333333333333327  5.881304775837604e-16  \\
            5.333333333333327  5.881304775837604e-16  \\
            5.333333333333327  5.881304775837604e-16  \\
            5.66666666666667  9.508154051376487e-16  \\
        }
        ;
    \node[]  at (axis cs:6.27,1.0e-8) {$r_{9,9}$};
    \node[]  at (axis cs:2.3333333333333366,0.0001774418739862507) {$r_{2,2}$};
    \node[]  at (axis cs:2.3333333333333366,0.0001774418739862507) {$r_{2,2}$};
    \node[]  at (axis cs:2.3333333333333366,0.0001774418739862507) {$r_{2,2}$};
    \node[]  at (axis cs:3.333333333333332,2.2834963464047978e-7) {$r_{3,3}$};
    \node[]  at (axis cs:3.333333333333332,2.2834963464047978e-7) {$r_{3,3}$};
    \node[]  at (axis cs:3.6666666666666696,1.6144428882142386e-10) {$r_{4,4}$};
    \node[]  at (axis cs:3.6666666666666696,1.6144428882142386e-10) {$r_{4,4}$};
    \node[]  at (axis cs:3.6666666666666696,1.6144428882142386e-10) {$r_{4,4}$};
    \node[]  at (axis cs:4.333333333333328,7.229460053600613e-14) {$r_{5,5}$};
    \node[]  at (axis cs:4.333333333333328,7.229460053600613e-14) {$r_{5,5}$};
    \node[]  at (axis cs:5.0,6.1485062293324904e-15) {$r_{6,6}$};
    \node[]  at (axis cs:5.0,6.1485062293324904e-15) {$r_{6,6}$};
    \node[]  at (axis cs:5.0,6.1485062293324904e-15) {$r_{6,6}$};
    \node[]  at (axis cs:5.333333333333327,2.940652387918802e-15) {$r_{7,7}$};
    \node[]  at (axis cs:5.333333333333327,2.940652387918802e-15) {$r_{7,7}$};
    \node[]  at (axis cs:5.333333333333327,2.940652387918802e-15) {$r_{7,7}$};
    \node[]  at (axis cs:5.66666666666667,4.754077025688244e-15) {$r_{8,8}$};
\end{axis}
\end{tikzpicture}%
	\end{subfigure}

	\begin{subfigure}[b]{0.48\textwidth}
		\raggedleft
		\begin{tikzpicture}[/tikz/background rectangle/.style={fill={rgb,1:red,1.0;green,1.0;blue,1.0}, fill opacity={1.0}, draw opacity={1.0}}, show background rectangle]
\begin{axis}[point meta max={nan}, point meta min={nan}, legend cell align={left}, legend columns={1}, title={$\lVert hA \rVert_{1}=10.0$}, title style={at={{(0.5,1)}}, anchor={south}, color={rgb,1:red,0.0;green,0.0;blue,0.0}, draw opacity={1.0}, rotate={0.0}, align={center}}, legend style={color={rgb,1:red,0.0;green,0.0;blue,0.0}, draw opacity={1.0}, line width={1}, solid, fill={rgb,1:red,1.0;green,1.0;blue,1.0}, fill opacity={1.0}, text opacity={1.0}, text={rgb,1:red,0.0;green,0.0;blue,0.0}, cells={anchor={center}}, at={(1.02, 1)}, anchor={north west}}, axis background/.style={fill={rgb,1:red,1.0;green,1.0;blue,1.0}, opacity={1.0}}, anchor={north west}, scaled x ticks={false}, xlabel={cost,\ [$C$]}, x tick style={color={rgb,1:red,0.0;green,0.0;blue,0.0}, opacity={1.0}}, x tick label style={color={rgb,1:red,0.0;green,0.0;blue,0.0}, opacity={1.0}, rotate={0}}, xmajorgrids={true}, xmin={4.0}, xmax={8.453333333333344}, xticklabels={{$5$,$6$,$7$,$8$}}, xtick={{5.0,6.0,7.0,8.0}}, xtick align={inside}, xticklabel style={color={rgb,1:red,0.0;green,0.0;blue,0.0}, draw opacity={1.0}, rotate={0.0}}, x grid style={color={rgb,1:red,0.0;green,0.0;blue,0.0}, draw opacity={0.1}, line width={0.5}, solid}, axis x line*={left}, x axis line style={color={rgb,1:red,0.0;green,0.0;blue,0.0}, draw opacity={1.0}, line width={1}, solid}, scaled y ticks={false}, ylabel={$\lVert w_{\alpha}(h A)-e^{h A} \rVert_{1} / (\lVert hA \rVert_{1} \cdot \lVert e^{hA} \rVert_{1})$}, y tick style={color={rgb,1:red,0.0;green,0.0;blue,0.0}, opacity={1.0}}, y tick label style={color={rgb,1:red,0.0;green,0.0;blue,0.0}, opacity={1.0}, rotate={0}}, ymode={log}, log basis y={10}, ymajorgrids={true}, ymin={1.0e-17}, ymax={1.0}, yticklabels={{$10^{-16}$,$10^{-8}$,$10^{-12}$,$10^{-4}$}}, ytick={{1.0e-16,1.0e-8,1.0e-12,0.0001}}, ytick align={inside}, yticklabel style={color={rgb,1:red,0.0;green,0.0;blue,0.0}, draw opacity={1.0}, rotate={0.0}}, y grid style={color={rgb,1:red,0.0;green,0.0;blue,0.0}, draw opacity={0.1}, line width={0.5}, solid}, axis y line*={left}, y axis line style={color={rgb,1:red,0.0;green,0.0;blue,0.0}, draw opacity={1.0}, line width={1}, solid}, colorbar={false}]
    \addplot[color={rgb,1:red,0.0;green,0.6056;blue,0.9787}, name path={5ce7d496-1ce1-4223-8ecf-924ffc18cf9d}, const plot, draw opacity={1.0}, line width={2}, dashed]
        table[row sep={\\}]
        {
            \\
            4.333333333333328  1.0  \\
            4.333333333333328  0.1  \\
            5.333333333333327  0.010000000000000002  \\
            5.66666666666667  0.001  \\
            5.66666666666667  0.0001  \\
            6.333333333333335  1.0e-5  \\
            6.666666666666664  1.0e-6  \\
            6.666666666666664  1.0e-7  \\
            7.333333333333339  1.0e-8  \\
            7.333333333333339  1.0e-9  \\
            7.333333333333339  1.0e-10  \\
            7.66666666666666  1.0e-11  \\
            7.66666666666666  1.0e-12  \\
            8.333333333333343  1.0e-13  \\
            8.333333333333343  1.0e-14  \\
            8.333333333333343  1.0e-15  \\
            8.333333333333343  1.0e-16  \\
        }
        ;
    \addplot[color={rgb,1:red,0.0;green,0.0;blue,0.0}, name path={512ed070-fcb6-4570-85cc-8fbaecd14a0a}, draw opacity={1.0}, line width={1}, dotted]
        table[row sep={\\}]
        {
            \\
            8.333333333333332  1.0e-34  \\
            8.333333333333332  1.0e17  \\
        }
        ;
    \addplot[color={rgb,1:red,0.2422;green,0.6433;blue,0.3044}, name path={7a37b040-48ff-43d6-afa2-7ef13ce63ab4}, const plot, draw opacity={1.0}, line width={2}, solid]
        table[row sep={\\}]
        {
            \\
            4.333333333333328  0.0007961380756841287  \\
            4.333333333333328  0.0007961380756841287  \\
            5.333333333333327  7.790786301733951e-8  \\
            5.66666666666667  2.5996321045477483e-8  \\
            5.66666666666667  2.5996321045477483e-8  \\
            6.333333333333335  3.1657809081515834e-12  \\
            6.666666666666664  1.3281386646090103e-14  \\
            6.666666666666664  1.3281386646090103e-14  \\
            7.333333333333339  3.557891647839391e-16  \\
            7.333333333333339  1.6262822718640128e-16  \\
            7.333333333333339  2.624630722740447e-16  \\
            7.66666666666666  3.5141447454905857e-16  \\
            7.66666666666666  3.5141447454905857e-16  \\
            8.333333333333343  1.5171183841618588e-16  \\
            8.333333333333343  1.5171183841618588e-16  \\
            8.333333333333343  1.5171183841618588e-16  \\
            8.333333333333343  1.5171183841618588e-16  \\
        }
        ;
    \node[]  at (axis cs:8.249999999999998,1.0e-8) {$r_{13,13}$};
    \node[]  at (axis cs:4.333333333333328,0.003980690378420644) {$r_{2,2}$};
    \node[]  at (axis cs:4.333333333333328,0.003980690378420644) {$r_{2,2}$};
    \node[]  at (axis cs:5.333333333333327,3.8953931508669756e-7) {$r_{5,5}$};
    \node[]  at (axis cs:5.66666666666667,1.2998160522738743e-7) {$r_{4,4}$};
    \node[]  at (axis cs:5.66666666666667,1.2998160522738743e-7) {$r_{4,4}$};
    \node[]  at (axis cs:6.333333333333335,1.582890454075792e-11) {$r_{7,7}$};
    \node[]  at (axis cs:6.666666666666664,6.640693323045051e-14) {$r_{8,8}$};
    \node[]  at (axis cs:6.666666666666664,6.640693323045051e-14) {$r_{8,8}$};
    \node[]  at (axis cs:7.333333333333339,8.7789458239196954e-12) {$r_{13,13}$};
    \node[]  at (axis cs:7.333333333333339,8.1314113593200645e-16) {$r_{9,9}$};
    \node[]  at (axis cs:7.333333333333339,8.3123153613702236e-14) {$r_{7,7}$};
    \node[]  at (axis cs:7.66666666666666,1.7570723727452929e-15) {$r_{8,8}$};
    \node[]  at (axis cs:7.66666666666666,1.7570723727452929e-15) {$r_{8,8}$};
    \node[]  at (axis cs:8.333333333333343,7.585591920809295e-16) {$r_{13,13}$};
    \node[]  at (axis cs:8.333333333333343,7.585591920809295e-16) {$r_{13,13}$};
    \node[]  at (axis cs:8.333333333333343,7.585591920809295e-16) {$r_{13,13}$};
    \node[]  at (axis cs:8.333333333333343,7.585591920809295e-16) {$r_{13,13}$};
\end{axis}
\end{tikzpicture}%
	\end{subfigure}
	~
	\begin{subfigure}[b]{0.48\textwidth}
		\raggedright
		\begin{tikzpicture}[/tikz/background rectangle/.style={fill={rgb,1:red,1.0;green,1.0;blue,1.0}, fill opacity={1.0}, draw opacity={1.0}}, show background rectangle]
\begin{axis}[point meta max={nan}, point meta min={nan}, legend cell align={left}, legend columns={1}, title={$\lVert hA \rVert_{1}=100.0$}, title style={at={{(0.5,1)}}, anchor={south}, color={rgb,1:red,0.0;green,0.0;blue,0.0}, draw opacity={1.0}, rotate={0.0}, align={center}}, legend style={color={rgb,1:red,0.0;green,0.0;blue,0.0}, draw opacity={1.0}, line width={1}, solid, fill={rgb,1:red,1.0;green,1.0;blue,1.0}, fill opacity={1.0}, text opacity={1.0}, text={rgb,1:red,0.0;green,0.0;blue,0.0}, cells={anchor={center}}, at={(1.02, 1)}, anchor={north west}}, axis background/.style={fill={rgb,1:red,1.0;green,1.0;blue,1.0}, opacity={1.0}}, anchor={north west}, scaled x ticks={false}, xlabel={cost,\ [$C$]}, x tick style={color={rgb,1:red,0.0;green,0.0;blue,0.0}, opacity={1.0}}, x tick label style={color={rgb,1:red,0.0;green,0.0;blue,0.0}, opacity={1.0}, rotate={0}}, xmajorgrids={true}, xmin={7.0}, xmax={12.483333333333333}, xticklabels={{$8$,$9$,$10$,$11$,$12$}}, xtick={{8.0,9.0,10.0,11.0,12.0}}, xtick align={inside}, xticklabel style={color={rgb,1:red,0.0;green,0.0;blue,0.0}, draw opacity={1.0}, rotate={0.0}}, x grid style={color={rgb,1:red,0.0;green,0.0;blue,0.0}, draw opacity={0.1}, line width={0.5}, solid}, axis x line*={left}, x axis line style={color={rgb,1:red,0.0;green,0.0;blue,0.0}, draw opacity={1.0}, line width={1}, solid}, scaled y ticks={false}, ylabel={$\lVert w_{\alpha}(h A)-e^{h A} \rVert_{1} / (\lVert hA \rVert_{1} \cdot \lVert e^{hA} \rVert_{1})$}, y tick style={color={rgb,1:red,0.0;green,0.0;blue,0.0}, opacity={1.0}}, y tick label style={color={rgb,1:red,0.0;green,0.0;blue,0.0}, opacity={1.0}, rotate={0}}, ymode={log}, log basis y={10}, ymajorgrids={true}, ymin={1.0e-17}, ymax={1.0}, yticklabels={{$10^{-16}$,$10^{-8}$,$10^{-12}$,$10^{-4}$}}, ytick={{1.0e-16,1.0e-8,1.0e-12,0.0001}}, ytick align={inside}, yticklabel style={color={rgb,1:red,0.0;green,0.0;blue,0.0}, draw opacity={1.0}, rotate={0.0}}, y grid style={color={rgb,1:red,0.0;green,0.0;blue,0.0}, draw opacity={0.1}, line width={0.5}, solid}, axis y line*={left}, y axis line style={color={rgb,1:red,0.0;green,0.0;blue,0.0}, draw opacity={1.0}, line width={1}, solid}, colorbar={false}]
    \addplot[color={rgb,1:red,0.0;green,0.6056;blue,0.9787}, name path={66c0daa0-b1b5-4f03-b94a-b4ca845da953}, const plot, draw opacity={1.0}, line width={2}, dashed]
        table[row sep={\\}]
        {
            \\
            7.333333333333339  1.0  \\
            8.333333333333343  0.1  \\
            8.333333333333343  0.010000000000000002  \\
            8.666666666666655  0.001  \\
            9.333333333333346  0.0001  \\
            9.666666666666652  1.0e-5  \\
            10.0  1.0e-6  \\
            10.333333333333348  1.0e-7  \\
            10.333333333333348  1.0e-8  \\
            10.666666666666654  1.0e-9  \\
            10.666666666666654  1.0e-10  \\
            11.33333333333334  1.0e-11  \\
            11.33333333333334  1.0e-12  \\
            11.33333333333334  1.0e-13  \\
            11.33333333333334  1.0e-14  \\
            11.666666666666663  1.0e-15  \\
            12.333333333333332  1.0e-16  \\
        }
        ;
    \addplot[color={rgb,1:red,0.0;green,0.0;blue,0.0}, name path={fd6c6cfa-b11b-4e54-97c0-616ef77f4abc}, draw opacity={1.0}, line width={1}, dotted]
        table[row sep={\\}]
        {
            \\
            12.333333333333332  1.0e-34  \\
            12.333333333333332  1.0e17  \\
        }
        ;
    \addplot[color={rgb,1:red,0.2422;green,0.6433;blue,0.3044}, name path={c5ce69e0-d3c2-486a-933c-7e0a8e9f867a}, const plot, draw opacity={1.0}, line width={2}, solid]
        table[row sep={\\}]
        {
            \\
            7.333333333333339  0.001600067433131433  \\
            8.333333333333343  6.409832779564424e-7  \\
            8.333333333333343  0.00010130468218501036  \\
            8.666666666666655  1.3248762921382483e-7  \\
            9.333333333333346  6.2528470006245e-11  \\
            9.666666666666652  4.089163139147262e-13  \\
            10.0  1.6274174523002082e-12  \\
            10.333333333333348  2.126987947587547e-15  \\
            10.333333333333348  3.4951591195387613e-15  \\
            10.666666666666654  2.0607087069943123e-16  \\
            10.666666666666654  2.0607087069943123e-16  \\
            11.33333333333334  9.820790169143554e-17  \\
            11.33333333333334  9.820790169143554e-17  \\
            11.33333333333334  9.820790169143554e-17  \\
            11.33333333333334  9.820790169143554e-17  \\
            11.666666666666663  2.749966880554826e-16  \\
            12.333333333333332  5.473696265830035e-17  \\
        }
        ;
    \node[]  at (axis cs:12.209999999999999,1.0e-8) {$r_{13,13}$};
    \node[]  at (axis cs:7.333333333333339,0.008000337165657165) {$r_{2,2}$};
    \node[]  at (axis cs:8.333333333333343,3.2049163897822117e-6) {$r_{5,5}$};
    \node[]  at (axis cs:8.333333333333343,0.0005065234109250518) {$r_{2,2}$};
    \node[]  at (axis cs:8.666666666666655,6.624381460691241e-7) {$r_{4,4}$};
    \node[]  at (axis cs:9.333333333333346,3.1264235003122503e-10) {$r_{7,7}$};
    \node[]  at (axis cs:9.666666666666652,2.0445815695736312e-12) {$r_{8,8}$};
    \node[]  at (axis cs:10.0,8.137087261501041e-12) {$r_{6,6}$};
    \node[]  at (axis cs:10.333333333333348,1.0634939737937735e-14) {$r_{9,9}$};
    \node[]  at (axis cs:10.333333333333348,1.7475795597693807e-14) {$r_{7,7}$};
    \node[]  at (axis cs:10.666666666666654,1.0303543534971562e-15) {$r_{8,8}$};
    \node[]  at (axis cs:10.666666666666654,1.0303543534971562e-15) {$r_{8,8}$};
    \node[]  at (axis cs:11.33333333333334,4.910395084571777e-16) {$r_{13,13}$};
    \node[]  at (axis cs:11.33333333333334,4.910395084571777e-16) {$r_{13,13}$};
    \node[]  at (axis cs:11.33333333333334,4.910395084571777e-16) {$r_{13,13}$};
    \node[]  at (axis cs:11.33333333333334,4.910395084571777e-16) {$r_{13,13}$};
    \node[]  at (axis cs:11.666666666666663,1.374983440277413e-15) {$r_{8,8}$};
    \node[]  at (axis cs:12.333333333333332,2.7368481329150173e-16) {$r_{13,13}$};
\end{axis}
\end{tikzpicture}%
	\end{subfigure}
	\caption{Comparison of tolerance (dashed) and practical average cost plotted vs. total computational cost of diagonal methods.
		The thin vertical dotted line represents the cost in the current software that implements \cite{higham09tsa}.
		Note that the distance between the dotted vertical line and the solid one represents cost savings.}
	\label{fig:error_vs_cost_diag}
\end{figure}

Then, for every norm--tolerance pair in $h=10^{m}, \ m=-3,-2,\ldots,2$ and $tol=10^{-k}, \ k=0,1,\ldots,16$ the algorithm returns $w_{\alpha}(h A)$ and the corresponding computational cost in terms of matrix--matrix products, \( C \).
For each norm \( h \), we calculate the normalized relative error
\begin{equation}\label{eq.norm_relerr}
	\frac{\opnorm{w_{\alpha}(h A)-\e^{h A}}}{\opnorm{A} \cdot \opnorm{\e^{h A}}},
\end{equation}
where the reference \( \e^{h A} \) is computed numerically in arbitrary arithmetic, obtaining error--cost pairs.
Then in \autoref{fig:error_vs_cost_superdiag}, we plot these pairs alongside with the theoretical estimates from \autoref{tab.thetaAll}.
The same numerical tests for the diagonal methods listed in \autoref{tab.thetaDiag} are presented in \autoref{fig:error_vs_cost_diag}.

In \autoref{fig:error_vs_cost_superdiag} one can see that the average error of superdiagonal methods from \autoref{tab.thetaAll} is below the tolerance predicted by the backward error analysis \eqref{bound.2}.
The only exception is $r_{12,8}(x)$	for low tolerance, where $r_{13,13}(x)$ should be preferred.

The results that one would obtain using the schemes from the current numerical software are included for comparison as vertical dotted lines, which correspond to unchanging approximation to the round-off error.
The specific scheme used by the algorithm is annotated above the error line.
More specifically, if we take $m=-1$, i.e. $\opnorm{A} \leq 0.1$, one of six methods seen in the \autoref{fig:error_vs_cost_superdiag} will be used, while $r_{5,5}$ would be used at a higher cost in $\textsc{MATLAB}$ or Julia.

\paragraph{Numerical example 2: Symplectic}
Let $J$ be the symplectic matrix and $B$ a symmetric matrix, then the Hamiltonian matrix $A=JB$ belongs to the symplectic algebra and
\[
\left( \e^A \right)^T J \e^A=J, \qquad
J = \begin{bmatrix}
0  & I  \\
-I & 0
\end{bmatrix}.
\]
This property is preserved with diagonal Pad\'e approximants, and we can analyse both the accuracy of the methods in the \autoref{tab.thetaDiag} versus computational cost and to check that the symplectic property is preserved.

\begin{figure}[p]
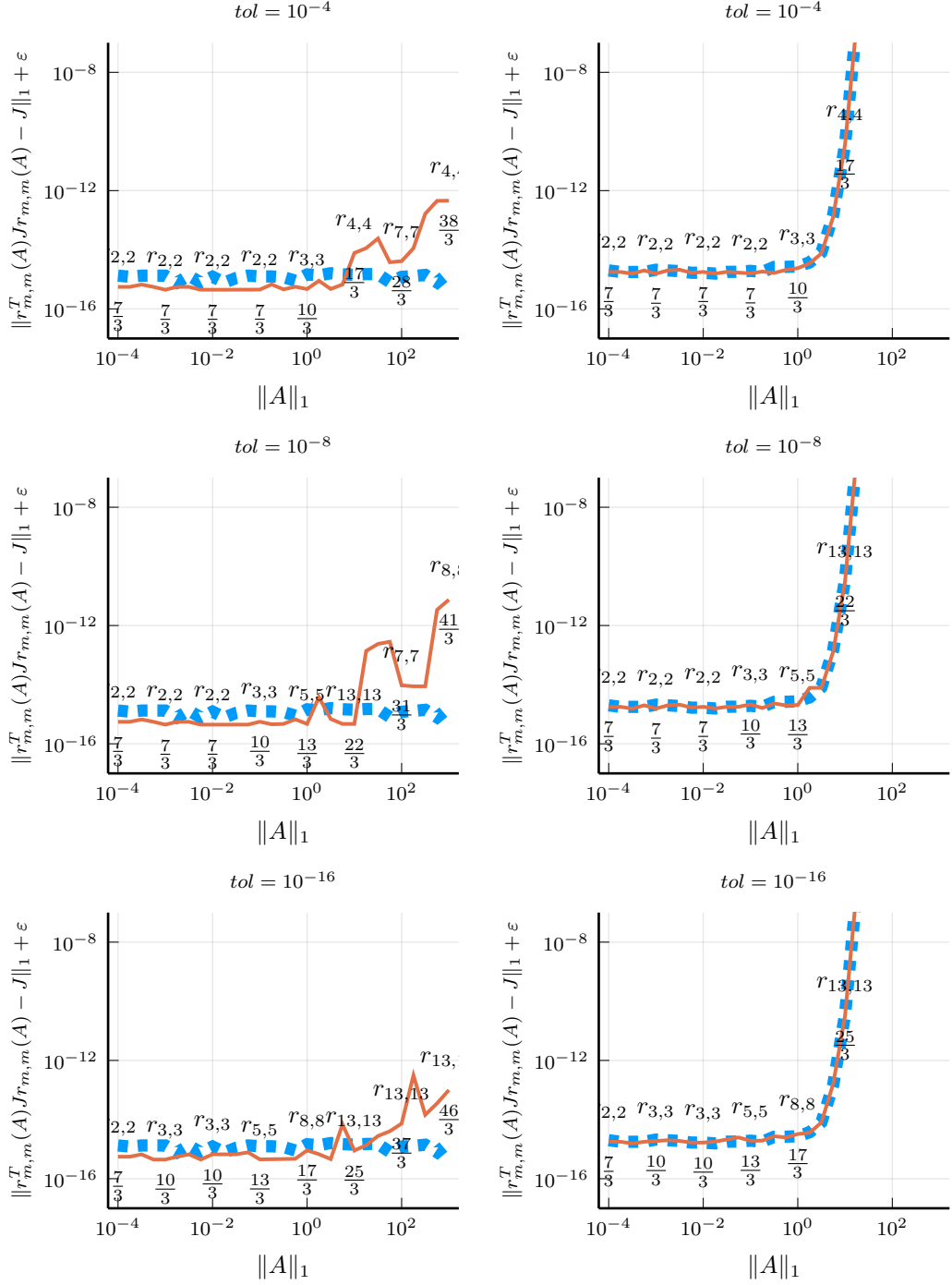

\begin{subfigure}[b]{0.48\textwidth}
	\raggedleft
	\input{fig/sympl_err_vs_cost_best_case.tikz}%
\end{subfigure}
~
\begin{subfigure}[b]{0.48\textwidth}
	\raggedright
	\input{fig/sympl_err_vs_cost_rand.tikz}%
\end{subfigure}
\caption{Comparison of symplectic property preservation by the algorithm using only diagonal Padé methods for matrices \eqref{eq:sympl_err_best_case} (left) and \eqref{eq:sympl_err_rand} (right).
	Reference (dashed) is the \texttt{LinearAlgebra.exp} function from Julia 1.10.2.
	Machine \( \varepsilon \) is added to avoid singularities.}
\label{fig:sympl_err}
\end{figure}

To illustrate the preservation of the symplectic property, we fix three tolerance values ($10^{-4}$, $10^{-8}$, $10^{-16}$) to calculate \( w_{\alpha}(A) \) and two types of matrices.
As a simple case we use the following matrix with diagonal blocks:
\begin{equation}
\label{eq:sympl_err_best_case}
A = \begin{bmatrix}
	0  & D  \\
	-D & 0
\end{bmatrix},
\quad D = diag(-26,-25,...,26) \in \mathbb{R}^{53 \times 53}.
\end{equation}
As a generic example, we take a Hamiltonian matrix \( A \) with elements sampled from \( \left[ -1;\,1 \right] \):
\begin{equation}
\label{eq:sympl_err_rand}
A = \begin{bmatrix}
	F & H       \\
	G & -F^{T}
\end{bmatrix},
\quad F,G,H \in \mathbb{R}^{53 \times 53},
\end{equation}
where \( G,H \) are symmetric.
As before, \( A \) is first normalized, and then we use the proposed algorithm to calculate the standard relative error w.r.t. $J$
\begin{equation}\label{eq.sympl_err}
	\frac{\opnorm{w_{\alpha}(A)^{T} J w_{\alpha}(A) - J}}{\opnorm{J}}
	= \opnorm{r_{m,m}^{T}(A) J r_{m,m}(A) - J},
\end{equation}
for a set of norms.
Then for each tolerance we plot the error \eqref{eq.sympl_err} versus the matrix norm indicating the cost and the corresponding method.
It can be seen that for each of the prescribed tolerances the algorithm preserves the  on par with the reference \( r_{13,13} \) while requiring less computational effort.
From both \autoref{fig:error_vs_cost_diag} and \autoref{fig:sympl_err} one can see that for all considered tolerance values there are computational savings while the symplectic property is preserved on par with \( r_{13,13} \).

We attribute the discrepancies for larger norms to a simpler implementation \( r_{13,13}(x) \) in our experimental code and Julia's \texttt{LinearAlgebra} which, e.g. uses matrix rebalancing and directly calls some \textsc{BLAS} functions.
A highly-optimized implementation could switch to the built-in methods when necessary.

\begin{figure}[p]
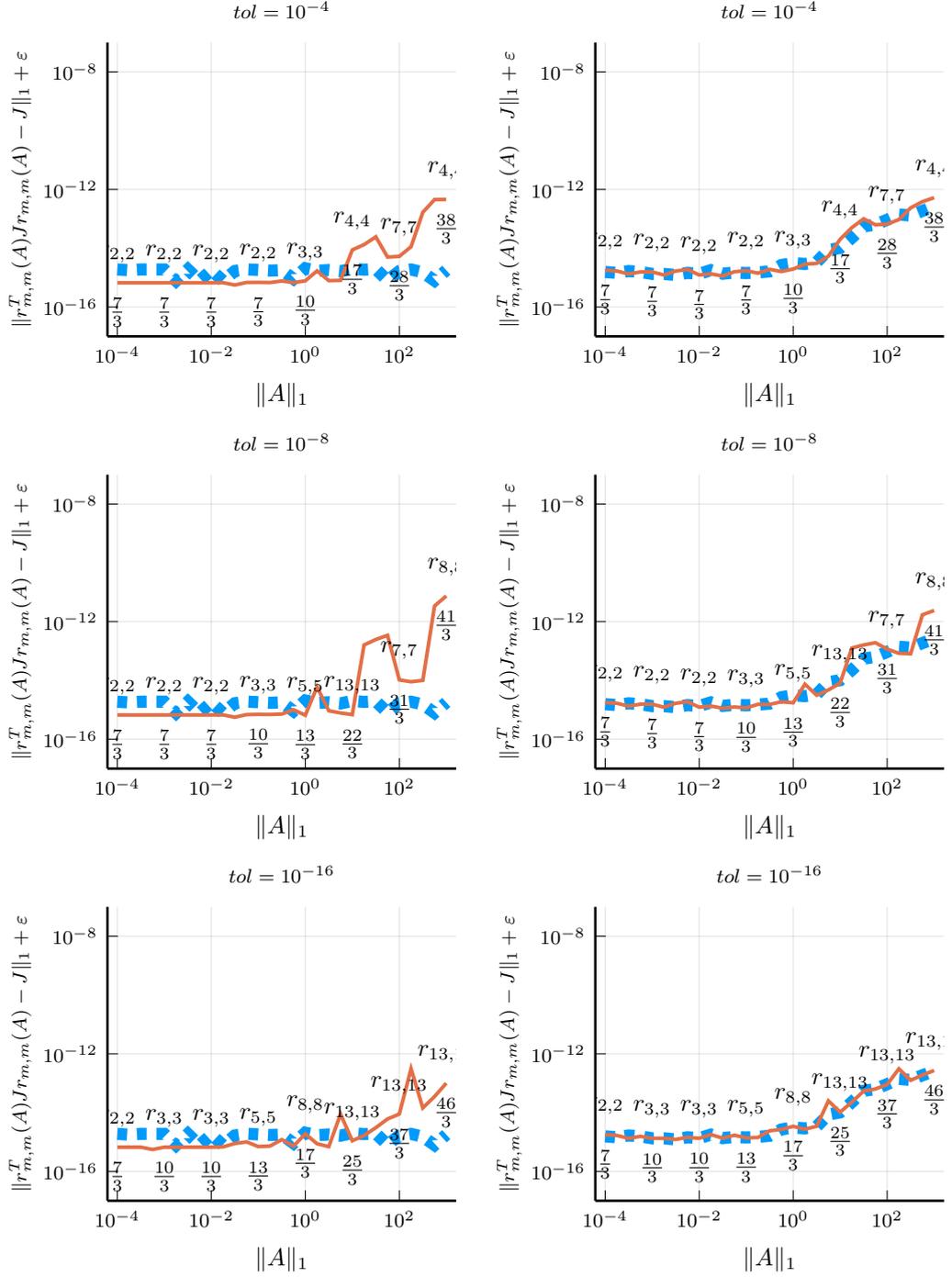

\begin{subfigure}[b]{0.48\textwidth}
	\raggedleft
	\input{fig/unit_err_vs_cost_best_case.tikz}%
\end{subfigure}
~
\begin{subfigure}[b]{0.48\textwidth}
	\raggedright
	\input{fig/unit_err_vs_cost_rand.tikz}%
\end{subfigure}
\caption{Comparison of unitarity preservation (where \( J=I \)) by the algorithm using only diagonal Padé methods for matrices \eqref{eq:sympl_err_best_case} (left) and a random complex matrix (right).
	Reference (dashed) is the \texttt{LinearAlgebra.exp} function from Julia 1.10.2.
	Machine \( \varepsilon \) is added to avoid singularities.}
\label{fig:unit_err}
\end{figure}

\paragraph{Numerical example 3: Unitary}
Similarly to the sypmlectic case, we consider the unitary one.
To do so, we use two matrices: the previous \eqref{eq:sympl_err_best_case} as well as $A=iB+C,\ A \in \mathbb{C}^{101 \times 101}$ with $B$ symmetric and $C$ skew-symmetric matrices.
We measure the relative error
\begin{equation}\label{eq.unit_err}
	\frac{\opnorm{w_{\alpha}(A)^{\dagger} w_{\alpha}(A) - I}}{\opnorm{I}}
	= \opnorm{r_{m,m}^{\dagger}(A) r_{m,m}(A) - I},
\end{equation}
and the results are represented in \autoref{fig:unit_err}.
As in the previous example, the property of interest is preserved, and there is  computational cost reduction.

For this problem, if the accuracy in the solution is also required to round-off, tailored polynomial approximations can be very efficient schemes \cite{bader22aea}.

\section{Concluding remarks}
\label{sec.7}
We have presented a new algorithm to approximate the matrix exponential function for different tolerances at a lower computational cost than the standard methods for a wide range of matrices.
The algorithm can be easily adjusted to use only Taylor methods (i.e., matrix--matrix products without inverses) for problems where products are considerably cheaper than the inverse or to use only diagonal Padé approximants to preserve the Lie group structure.

In some cases where $\opnorm{A^k}^{1/k} \ll \opnorm{A}, \ k > 1$ one can extend the error analysis to reduce the cost by allowing to use cheaper schemes or to reduce the number of scalings.
Since our algorithm computes several powers of $A$, this can be easily checked and, if the user is interested in this class of problems, it can be easily adjusted for them following \cite{higham09tsa,preprint,bader19ctm}.
We suggest leaving it only as an option otherwise it adds some extra cost for all problems, and in many cases this is not necessary.

The analysis presented in this paper can also be applied to the approximation of other matrix functions \cite{higham08fom} where, in general, diagonal Pad\'e approximants are no longer symmetric, i.e. $q_{m,m}(x)\neq p_{m,m}(-x)$.
Obviously, if the matrix has some particular structure like a perturbation of a matrix whose exponential is easy and cheap to compute, other techniques should be used \cite{bader15tss}.

We have decomposed several Pad\'e approximants into fractions of lower degree that can be computed independently, making them suitable for parallelization. We will carry a deep study of the computation of the matrix exponential to be computed in parallel following \cite{blanes22preprint}.

\subsection*{Acknowledgements}
This work has been funded by Ministerio de Ciencia e Innovacion (Spain) through project PID2022-136585NB-C21, MCIN/AEI/10.13039/501100011033/FEDER, UE, and also by Generalitat Valenciana (Spain) through project CIAICO/2021/180.


\begin{thebibliography}{99}
\bibitem{almohy09ans} {A.H. Al-Mohy and N.J. Higham},
A new scaling and squaring algorithm for the matrix exponential,
SIAM J. Matrix Anal. Appl., \textbf{31}, (2009), pp. 970--989.

\bibitem{bader22aea} P. Bader, S. Blanes, F. Casas, M. Seydaoğlu,
An efficient algorithm to compute the exponential of skew-Hermitian matrices for the time integration of the Schrödinger equation,
Math. Comput. Sim., \textbf{194}, (2022), pp. 383--400.

\bibitem{bader15tss} {P. Bader, S. Blanes, and M. Seydao\u{g}lu},
The scaling, splitting and squaring method for the exponential of perturbed matrices,
SIAM  J. Matrix Anal. Appl., \textbf{36}, (2015), pp. 594--614.

\bibitem{preprint} {P. Bader, S. Blanes, and F. Casas},
An improved algorithm to compute the exponential of a matrix, arXiv:1710.10989 [math.NA], (2017), preprint.

\bibitem{bader19ctm} {P. Bader, S. Blanes, and F. Casas},
Computing the matrix exponential with an optimized Taylor polynomial approximation,
Mathematics, \textbf{7} (2019), 1174.

\bibitem{julia17} J. Bezanson, A. Edelman, S. Karpinski, \& V. Shah, Julia: A fresh approach to numerical computing,
{SIAM Review}. \textbf{59}, 65-98 (2017).

\bibitem{blanes22preprint} {S. Blanes},
Parallel computation of functions of matrices and their action on vectors
arXiv:2210.03714 [math.NA], (2022), preprint.

\bibitem{blanes09tme}
S.~Blanes, F.~Casas, J.A. Oteo, and J.~Ros,
\newblock The {M}agnus expansion and some of its applications.
\newblock Phys. Rep., \textbf{470} (2009), pp. 151--238.

\bibitem{celledoni00ate} {E. Celledoni and A. Iserles},
{Approximating the exponential from a {L}ie  algebra to a {L}ie group},
Math. Comput., \textbf{69} (2000), pp.  {1457--1480}.

\bibitem{celledoni01mft}  {E. Celledoni and A. Iserles},
{Methods for the approximation of the matrix exponential in a {L}ie-algebraic setting},
IMA J. Numer. Anal.,  \textbf{21} (2001), pp. {463--488}.

\bibitem{cody69cra} W.J. Cody, G. Meinardus, and R. S. Varga,
Chebyshev rational approximation to $\e^{-x}$ in $[0,+\infty)$ and applications to heat-conduction problems,
J. Approximation Theory, 2 (1969), pp. 50--65.

\bibitem{crouch93nio} P.E. Crouch and R. Grossman, Numerical integration of ordinary differential equations on manifolds,
J. Nonlinear Sci. 3 (1993), 1–33.

\bibitem{higham05tsa} {N.J. Higham},
The scaling and squaring method for the matrix exponential revisited,
SIAM  J. Matrix Anal. Appl., \textbf{26}, (2005), pp. 1179--1193.

\bibitem{higham09tsa} {N.J. Higham},
The scaling and squaring method for the matrix exponential revisited,
SIAM Review \textbf{51} (2009), pp. 747-764.

\bibitem{higham08fom} {N.J. Higham},
{Functions of Matrices: Theory and Computation},
Society for Industrial and Applied Mathematics, Philadelphia, PA, USA (2008).

\bibitem{higham10cma} {N.J. Higham and A.H. Al-Mohy},
Computing matrix functions,
Acta Numerica, \textbf{51}, (2010), pp. 159--208.

\bibitem{hochbruck10ein} M. Hochbruck and A. Ostermann, Exponential integrators,
{Acta Numerica}, \textbf{19}, (2010),  pp. 209--286.

\bibitem{Iserles1982} A. Iserles, Composite exponential approximations. Mathematics Of Computation. \textbf{38}, 99-112 (1982), http://dx.doi.org/10.1090/S0025-5718-1982-0637289-7

\bibitem{iserles00lgm} {A. Iserles, H. Z. Munthe-Kaas, S.P. N{\o}rsett and A. Zanna}, Lie-group methods, {Acta Numerica}, \textbf{9}, (2000),  pp. 215--365.

\bibitem{lezcanocasado20aam} M. Lezcano-Casado,
Adaptive and momentum methods on manifolds through trivializations,
arXiv:2010.04617, (2020).

\bibitem{moler03ndw} {C.B. Moler and C.F. Van Loan},
Nineteen dubious ways to compute the exponential of a matrix, twenty-five years later,
SIAM Review, \textbf{45} (2003), pp. {3--49}.

\bibitem{munthekaas98rkm} H. Munthe-Kaas,
Runge--Kutta methods on Lie groups,
BIT 38 (1998), 92--111.

\bibitem{paterson73otn} M.S. Paterson and L.J. Stockmeyer,
On the number of nonscalar multiplications necessary to evaluate polynomials,
SIAM J. Comput., \textbf{2} (1973), pp. 60--66.

\bibitem{Rackauckas2017}  C. Rackauckas, Q. Nie, DifferentialEquations.jl – A Performant and Feature-Rich Ecosystem for Solving Differential Equations in Julia. The Journal Of Open Research Software. \textbf{5} (2017), https://doi.org/10.5334/jors.151

\bibitem{sastre12emr} {J. Sastre},
Efficient mixed rational and polynomial approximations of matrix functions,
Appl. Math. Comput., \textbf{218} (2012), pp. 11938--11946.

\bibitem{sastre18eeo} {J. Sastre},
Efficient evaluation of matrix polynomials,
Linear Algebra Appl., \textbf{539}, (2018), pp. 229--250.

\bibitem{sastre15nsq} {J. Sastre, J. Ib\'a\~nez, P. Ruiz, and E. Defez},
New scaling-squaring Taylor algorithms for computing the matrix exponential,
SIAM J. Sci. Comp., \textbf{37} (2015), pp. A439--A455.

\bibitem{sastre19btc} \blue{{J. Sastre, J. Ib\'a\~nez, and E. Defez},
	Boosting the computation of the matrix exponential,
	Appl. Math. Comput., \textbf{340} (2019), pp. 206-220.}

\bibitem{Ta2015} C. Ta, D. Wang, Q. Nie, An integration factor method for stochastic and stiff reaction–diffusion systems. Journal Of Computational Physics. \textbf{295} pp. 505-522 (2015,8), http://dx.doi.org/10.1016/j.jcp.2015.04.028

\bibitem{trefethen06tqa} L.N. Trefethen, J.A.C. Weideman, and T. Schmelzer,
Talbot quadratures and rational approximations,
BIT, 46 (2006), pp. 653--670.
\end{thebibliography}
\end{document}